\documentclass[final,3p]{elsarticle}
 \usepackage{graphics}
 \usepackage{graphicx}
 \usepackage{epsfig}
\usepackage{amssymb}
 \usepackage{amsthm}
 \usepackage{lineno}
 \usepackage{amsmath}
   \numberwithin{equation}{section}
\usepackage{mathrsfs}

\journal{Springer} 

\newtheorem{thm}{Theorem}[section]

\newtheorem{lem}[thm]{Lemma}

\biboptions{sort&compress,square}
\allowdisplaybreaks
\begin{document}
\begin{frontmatter}
\author{Hongfeng Li}
\ead{lihf728@nenu.edu.cn}
\author{Yong Wang\corref{cor2}}
\ead{wangy581@nenu.edu.cn}
\cortext[cor2]{Corresponding author.}

\address{School of Mathematics and Statistics, Northeast Normal University,
Changchun, 130024, China}

\title{The general Dabrowski-Sitarz-Zalecki type theorem for odd dimensional manifolds with boundary \uppercase\expandafter{\romannumeral2}}
\begin{abstract}
In \cite{WWW}, a general Dabrowski-Sitarz-Zalecki type theorem for odd dimensional manifolds with boundary was proved. In this paper, we give the proof of the another general Dabrowski-Sitarz-Zalecki type theorem for the spectral Einstein functional associated with the Dirac operator on odd dimensional manifolds with boundary.
\end{abstract}
\begin{keyword}Dirac operators; the spectral Einstein functional; the Dabrowski-Sitarz-Zalecki type theorem\\

\end{keyword}
\end{frontmatter}

\section{Introduction}
The noncommutative residue found in \cite{Gu,Wo} plays a prominent role in noncommutative geometry. For one-dimensional manifolds, the noncommutative residue was discovered by Adler \cite{MA} in connection with geometric aspects of nonlinear partial differential equations. For arbitrary closed compact $n$-dimensional manifolds, the noncommutative residue was introduced by Wodzicki in \cite{Wo} using the theory of zeta functions of elliptic pseudo-differential operators. In \cite{Co1}, Connes used the noncommutative residue to derive a conformal 4-dimensional Polyakov action analogy. Connes showed us that the noncommutative residue on a compact manifold $M$ coincided with the Dixmier's trace on pseudo-differential operators of order $-{\rm {dim}}M$ in \cite{Co2}.
And Connes claimed that the noncommutative residue of the square of the inverse of the Dirac operator was proportioned to the Einstein-Hilbert action.  Kastler \cite{Ka} gave a
brute-force proof of this theorem. Kalau and Walze proved this theorem in the normal coordinates system simultaneously in \cite{KW} .
Ackermann proved that
the Wodzicki residue  of the square of the inverse of the Dirac operator ${\rm  Wres}(D^{-2})$ in turn is essentially the second coefficient
of the heat kernel expansion of $D^{2}$ in \cite{Ac}.

On the other hand, Wang generalized the Connes' results to the case of manifolds with boundary in \cite{Wa1,Wa2},
and proved the Kastler-Kalau-Walze type theorem for the Dirac operator and the signature operator on lower-dimensional manifolds
with boundary \cite{Wa3}. In \cite{Wa3,Wa4}, Wang computed $\widetilde{{\rm Wres}}[\pi^+D^{-1}\circ\pi^+D^{-1}]$ and $\widetilde{{\rm Wres}}[\pi^+D^{-2}\circ\pi^+D^{-2}]$, where the two operators are symmetric, in these cases the boundary term vanished. But for $\widetilde{{\rm Wres}}[\pi^+D^{-1}\circ\pi^+D^{-3}]$, the authors got a nonvanishing boundary term \cite{Wa5}, and give a theoretical explanation for gravitational action on boundary. In \cite{DL}, the authors defined bilinear functionals of vector fields and differential forms, the densities of which yielded the  metric and Einstein tensors on even dimensional Riemannian manifolds. The authors computed the generalized noncommutative residue $\widetilde{{\rm Wres}}[\pi^+(c(X)D^{-1})\circ\pi^+(D^{-(2m-2)})]$, $\widetilde{{\rm Wres}}[\pi^+(\nabla_X^{S(TM)}D^{-1})\circ\pi^+(D^{-(2m-1)})]$ and $\widetilde{{\rm Wres}}[\pi^+(\nabla_X^{S(TM)}D^{-2})\circ\pi^+(D^{-(2m-2)})]$ on odd dimensional manifolds with boundary \cite{WWW}. In \cite{Wu1}, the authors defined the spectral Einstein functional associated with the Dirac operator for manifolds with boundary, and computed the noncommutative residue $\widetilde{{\rm Wres}}[\pi^+(\nabla_X^{S(TM)}\nabla_Y^{S(TM)}D^{-2})\circ\pi^+(D^{-(n-2)})]$  and $\widetilde{{\rm Wres}}[\pi^+(\nabla_X^{S(TM)}\nabla_Y^{S(TM)}D^{-1})\circ\pi^+(D^{-(n-1)})]$ on even and odd dimensional compact manifolds. {\bf The motivation of this paper} is to prove the another general Dabrowski-Sitarz-Zalecki type theorem associated with the Dirac operator for odd dimensional manifolds with boundary. That is, we want to compute ${\rm Wres}[\pi^+P_1\circ\pi^+P_2]$, where orders of $P_1,P_2$ are $a_1,a_2$ and $-a_1-a_2+2=n$ for odd dimensional manifolds with boundary. Motivated by \cite{Wa6,Wa7,Wu1}, we compute the noncommutative residue $\widetilde{{\rm Wres}}[\pi^+(c(Z)\nabla_X^{S(TM)}\nabla_Y^{S(TM)}D^{-1})\circ\pi^+(D^{-2m})]$  and $\widetilde{{\rm Wres}}[\pi^+(c(Z)\nabla_X^{S(TM)}\nabla_Y^{S(TM)}D^{-2})\circ\pi^+(D^{-(2m-1)})]$ on odd dimensional compact manifolds. Our main theorems are as follows.
\\
\begin{thm}\label{thmb1}
Let $M$ be an $n=(2m+1)$-dimensional oriented
compact spin manifold with the boundary $\partial M$, then we get the following equality:
\begin{align}
\label{b2773}
&\widetilde{{\rm Wres}}[\pi^+(c(Z)\nabla_X^{S(TM)}\nabla_Y^{S(TM)}D^{-1})\circ\pi^+(D^{-2m})]\nonumber\\
=&\int_{\partial M}\bigg\{\partial_{x_n}\Big(g(X^T,Y^T)Z_n\Big)\frac{\pi}{2(m+2)!}A_0-\partial_{x_n}\Big((X_nY_n)Z_n\Big)\frac{\pi}{(m+2)!}A_1+g(X^T,Y^T)Z_nh'(0)\Big(\frac{\pi i}{4(m+3)!}A_2\nonumber\\
&+\frac{\pi}{2(m+3)!}B_0+\frac{\pi}{4(m+2)!}C_0+\frac{\pi}{2(m+3)!}C_1-\frac{\pi}{2(m+3)!}D_1\Big)-X_nY_nZ_nh'(0)\Big(-\frac{\pi i}{2(m+3)!}A_3\nonumber\\
&+\frac{\pi}{(m+3)!}B_1+\frac{\pi}{2(m+2)!}C_2+\frac{\pi}{(m+3)!}C_3+\frac{\pi}{(m+3)!}D_2+\frac{\pi}{(m+2)!}D_3\Big)\nonumber\\
&-X(Y_n)Z_n\frac{\pi}{(m+1)!}D_0\bigg\}Vol(S^{n-2})2^md{\rm Vol_{M}},
\end{align}
where $X=\sum_{j=1}^nX_j\partial_{x_j}, Y=\sum_{l=1}^nY_l\partial_{x_l}, Z=\sum_{k=1}^nY_k\partial_{x_k}$ are vector fields, and $A_0$-$D_3$ are defined in section 3.
\end{thm}
\begin{thm}\label{cthmb1}
Let $M$ be an $n=(2m+1)$-dimensional oriented
compact spin manifold with the boundary $\partial M$, then we get the following equality:
\begin{align}
\label{cb2773}
&\widetilde{{\rm Wres}}[\pi^+(c(Z)\nabla_X^{S(TM)}\nabla_Y^{S(TM)}D^{-2})\circ\pi^+(D^{-(2m-1)})]\nonumber\\
=&\int_{\partial M}\bigg\{\partial_{x_n}\Big(g(X^T,Y^T)Z_n\Big)\frac{-\pi i}{2(m+2)!}E_0+\partial_{x_n}\Big((X_nY_n)Z_n\Big)\frac{\pi i}{(m+2)!}E_1\nonumber\\
&+g(X^T,Y^T)Z_nh'(0)\Big(\frac{\pi i}{4(m+3)!}E_2-\frac{\pi}{2(m+3)!}F_0+\frac{\pi}{8(m+3)!}G_0+\frac{3\pi}{16(m+2)!}H_1+\frac{\pi}{8(m+3)!}H_2\Big)\nonumber\\
&+X_nY_nZ_nh'(0)\Big(\frac{\pi i}{2(m+3)!}E_3+\frac{\pi i}{(m+3)!}F_1-\frac{\pi}{4(m+3)!}G_1+\frac{\pi}{8(m+2)!}H_3+\frac{\pi}{4(m+3)!}H_4\Big)\nonumber\\
&+X(Y_n)Z_n\frac{\pi}{(m+1)!}H_0\bigg\}Vol(S^{n-2})2^md{\rm Vol_{M}}.
\end{align}
where $X=\sum_{j=1}^nX_j\partial_{x_j}, Y=\sum_{l=1}^nY_l\partial_{x_l}, Z=\sum_{k=1}^nY_k\partial_{x_k}$ are vector fields, and $E_0$-$H_4$ are defined in section 4.
\end{thm}
\indent The paper is organized in the following way. In Sec.\ref{section:2}, we recall some basic facts and formulas about Boutet de
Monvel's calculus and the definition of the noncommutative residue for manifolds with boundary. In Sec.\ref{section:3} and Sec.\ref{section:4}, we define the spectral Einstein functional associated with the Dirac operator and prove the another Dabrowski-Sitarz-Zalecki type theorem for the spectral Einstein functional associated with the Dirac operator on odd dimensional manifolds with boundary.

\section{Boutet de Monvel's calculus}
\label{section:2}
 In this section, we recall some basic facts and formulas about Boutet de
Monvel's calculus and the definition of the noncommutative residue for manifolds with boundary which will be used in the following. For more details, see Section 2 in \cite{Wa3}.\\
 \indent Let $M$ be an n-dimensional compact oriented manifold with the boundary $\partial M$.
We assume that the metric $g^{TM}$ on $M$ has the following form near the boundary,
\begin{equation}
\label{b1}
g^{M}=\frac{1}{h(x_{n})}g^{\partial M}+dx _{n}^{2},
\end{equation}
where $g^{\partial M}$ is the metric on $\partial M$ and $h(x_n)\in C^{\infty}([0, 1)):=\{\widehat{h}|_{[0,1)}|\widehat{h}\in C^{\infty}((-\varepsilon,1))\}$ for
some $\varepsilon>0$ and $h(x_n)$ satisfies $h(x_n)>0$, $h(0)=1$, where $x_n$ denotes the normal directional coordinate. Let $U\subset M$ be a collar neighborhood of $\partial M$ which is diffeomorphic with $\partial M\times [0,1)$. By the definition of $h(x_n)\in C^{\infty}([0,1))$
and $h(x_n)>0$, there exists $\widehat{h}\in C^{\infty}((-\varepsilon,1))$ such that $\widehat{h}|_{[0,1)}=h$ and $\widehat{h}>0$ for some
sufficiently small $\varepsilon>0$. Then there exists a metric $g'$ on $\widetilde{M}=M\bigcup_{\partial M}\partial M\times
(-\varepsilon,0]$ which has the form on $U\bigcup_{\partial M}\partial M\times (-\varepsilon,0 ]$
\begin{equation}
\label{b2}
g'=\frac{1}{\widehat{h}(x_{n})}g^{\partial M}+dx _{n}^{2} ,
\end{equation}
such that $g'|_{M}=g$. We fix a metric $g'$ on the $\widetilde{M}$ such that $g'|_{M}=g$.

Let the Fourier transformation $F'$  be
\begin{equation*}
F':L^2({\bf R}_t)\rightarrow L^2({\bf R}_v);~F'(u)(v)=\int_\mathbb{R} e^{-ivt}u(t)dt
\end{equation*}
and let
\begin{equation*}
r^{+}:C^\infty ({\bf R})\rightarrow C^\infty (\widetilde{{\bf R}^+});~ f\rightarrow f|\widetilde{{\bf R}^+};~
\widetilde{{\bf R}^+}=\{x\geq0;x\in {\bf R}\}.
\end{equation*}
\indent We define $H^+=F'(\Phi(\widetilde{{\bf R}^+}));~ H^-_0=F'(\Phi(\widetilde{{\bf R}^-}))$ which satisfies
$H^+\bot H^-_0$, where $\Phi(\widetilde{{\bf R}^+}) =r^+\Phi({\bf R})$, $\Phi(\widetilde{{\bf R}^-}) =r^-\Phi({\bf R})$ and $\Phi({\bf R})$
denotes the Schwartz space. We have the following
 property: $h\in H^+~$ (resp. $H^-_0$) if and only if $h\in C^\infty({\bf R})$ which has an analytic extension to the lower (resp. upper) complex
half-plane $\{{\rm Im}\xi<0\}$ (resp. $\{{\rm Im}\xi>0\})$ such that for all nonnegative integer $l$,
 \begin{equation*}
\frac{d^{l}h}{d\xi^l}(\xi)\sim\sum^{\infty}_{k=1}\frac{d^l}{d\xi^l}(\frac{c_k}{\xi^k}),
\end{equation*}
as $|\xi|\rightarrow +\infty,{\rm Im}\xi\leq0$ (resp. ${\rm Im}\xi\geq0)$ and where $c_k\in\mathbb{C}$ are some constants.\\
 \indent Let $H'$ be the space of all polynomials and $H^-=H^-_0\bigoplus H';~H=H^+\bigoplus H^-.$ Denote by $\pi^+$ (resp. $\pi^-$) the
 projection on $H^+$ (resp. $H^-$). Let $\widetilde H=\{$rational functions having no poles on the real axis$\}$. Then on $\widetilde{H}$,
 \begin{equation}
 \label{b3}
\pi^+h(\xi_0)=\frac{1}{2\pi i}\lim_{u\rightarrow 0^{-}}\int_{\Gamma^+}\frac{h(\xi)}{\xi_0+iu-\xi}d\xi,
\end{equation}
where $\Gamma^+$ is a Jordan closed curve
included ${\rm Im}(\xi)>0$ surrounding all the singularities of $h$ in the upper half-plane and
$\xi_0\in {\bf R}$. In our computations, we only compute $\pi^+h$ for $h$ in $\widetilde{H}$. Similarly, define $\pi'$ on $\widetilde{H}$,
\begin{equation}
\label{b4}
\pi'h=\frac{1}{2\pi}\int_{\Gamma^+}h(\xi)d\xi.
\end{equation}
So $\pi'(H^-)=0$. For $h\in H\bigcap L^1({\bf R})$, $\pi'h=\frac{1}{2\pi}\int_{{\bf R}}h(v)dv$ and for $h\in H^+\bigcap L^1({\bf R})$, $\pi'h=0$.\\
\indent An operator of order $m\in {\bf Z}$ and type $d$ is a matrix\\
$$\widetilde{A}=\left(\begin{array}{lcr}
  \pi^+P+G  & K  \\
   T  &  \widetilde{S}
\end{array}\right):
\begin{array}{cc}
\   C^{\infty}(M,E_1)\\
 \   \bigoplus\\
 \   C^{\infty}(\partial{M},F_1)
\end{array}
\longrightarrow
\begin{array}{cc}
\   C^{\infty}(M,E_2)\\
\   \bigoplus\\
 \   C^{\infty}(\partial{M},F_2)
\end{array},
$$
where $M$ is a manifold with boundary $\partial M$ and
$E_1,E_2$~ (resp. $F_1,F_2$) are vector bundles over $M~$ (resp. $\partial M
$).~Here,~$P:C^{\infty}_0(\Omega,\overline {E_1})\rightarrow
C^{\infty}(\Omega,\overline {E_2})$ is a classical
pseudo-differential operator of order $m$ on $\Omega$, where
$\Omega$ is a collar neighborhood of $M$ and
$\overline{E_i}|M=E_i~(i=1,2)$. $P$ has an extension:
$~{\cal{E'}}(\Omega,\overline {E_1})\rightarrow
{\cal{D'}}(\Omega,\overline {E_2})$, where
${\cal{E'}}(\Omega,\overline {E_1})~({\cal{D'}}(\Omega,\overline
{E_2}))$ is the dual space of $C^{\infty}(\Omega,\overline
{E_1})~(C^{\infty}_0(\Omega,\overline {E_2}))$. Let
$e^+:C^{\infty}(M,{E_1})\rightarrow{\cal{E'}}(\Omega,\overline
{E_1})$ denote extension by zero from $M$ to $\Omega$ and
$r^+:{\cal{D'}}(\Omega,\overline{E_2})\rightarrow
{\cal{D'}}(\Omega, {E_2})$ denote the restriction from $\Omega$ to
$X$, then define
$$\pi^+P=r^+Pe^+:C^{\infty}(M,{E_1})\rightarrow {\cal{D'}}(\Omega,
{E_2}).$$

In addition, $P$ is supposed to have the
transmission property; this means that, for all $j,k,\alpha$, the
homogeneous component $p_j$ of order $j$ in the asymptotic
expansion of the
symbol $p$ of $P$ in local coordinates near the boundary satisfies:\\
$$\partial^k_{x_n}\partial^\alpha_{\xi'}p_j(x',0,0,+1)=
(-1)^{j-|\alpha|}\partial^k_{x_n}\partial^\alpha_{\xi'}p_j(x',0,0,-1),$$
then $\pi^+P:C^{\infty}(M,{E_1})\rightarrow C^{\infty}(M,{E_2})$. Let $G$, $T$ be respectively the singular Green operator
and the trace operator of order $m$ and type $d$. Let $K$ be a
potential operator and $S$ be a classical pseudo-differential
operator of order $m$ along the boundary. Denote by $B^{m,d}$ the collection of all operators of
order $m$
and type $d$,  and $\mathcal{B}$ is the union over all $m$ and $d$.\\
\indent Recall that $B^{m,d}$ is a Fr\'{e}chet space. The composition
of the above operator matrices yields a continuous map:
$B^{m,d}\times B^{m',d'}\rightarrow B^{m+m',{\rm max}\{
m'+d,d'\}}.$ Write $$\widetilde{A}=\left(\begin{array}{lcr}
 \pi^+P+G  & K \\
 T  &  \widetilde{S}
\end{array}\right)
\in B^{m,d},
 \widetilde{A}'=\left(\begin{array}{lcr}
\pi^+P'+G'  & K'  \\
 T'  &  \widetilde{S}'
\end{array} \right)
\in B^{m',d'}.$$\\
 The composition $\widetilde{A}\widetilde{A}'$ is obtained by
multiplication of the matrices (For more details see \cite{Sc}). For
example $\pi^+P\circ G'$ and $G\circ G'$ are singular Green
operators of type $d'$ and
$$\pi^+P\circ\pi^+P'=\pi^+(PP')+L(P,P').$$
Here $PP'$ is the usual
composition of pseudo-differential operators and $L(P,P')$ called
leftover term is a singular Green operator of type $m'+d$. For our case, $P,P'$ are classical pseudo-differential operators, in other words $\pi^+P\in \mathcal{B}^{\infty}$ and $\pi^+P'\in \mathcal{B}^{\infty}$ .\\
\indent Let $M$ be an $n$-dimensional compact oriented manifold with the boundary $\partial M$.
Denote by $\mathcal{B}$ the Boutet de Monvel's algebra. We recall that the main theorem in \cite{FGLS}.
\begin{thm}\label{th:32}{\rm\cite{FGLS}}{\bf(Fedosov-Golse-Leichtnam-Schrohe)}
 Let $X$ and $\partial X$ be connected, ${\rm dim}X=n\geq3$, and let $\widetilde{S}$ (resp. $\widetilde{S}'$) be the unit sphere about $\xi$ (resp. $\xi'$) and $\sigma(\xi)$ (resp. $\sigma(\xi')$) be the corresponding canonical
$n-1$ (resp. $(n-2)$) volume form.
 Set $\widetilde{A}=\left(\begin{array}{lcr}\pi^+P+G &   K \\
T &  \widetilde{S}    \end{array}\right)$ $\in \mathcal{B}$ , and denote by $p$, $b$ and $s$ the local symbols of $P, G$ and $\widetilde{S}$ respectively.
 Define:
 \begin{align}
{\rm{\widetilde{Wres}}}(\widetilde{A})=&\int_X\int_{\bf \widetilde{ S}}{\rm{tr}}_E\left[p_{-n}(x,\xi)\right]\sigma(\xi)dx \nonumber\\
&+2\pi\int_ {\partial X}\int_{\bf \widetilde{S}'}\left\{{\rm tr}_E\left[({\rm{tr}}b_{-n})(x',\xi')\right]+{\rm{tr}}
_F\left[s_{1-n}(x',\xi')\right]\right\}\sigma(\xi')dx',
\end{align}
where ${\rm{\widetilde{Wres}}}$ denotes the noncommutative residue of an operator in the Boutet de Monvel's algebra.\\
Then~~ a) ${\rm \widetilde{Wres}}([\widetilde{A},B])=0 $, for any
$\widetilde{A},B\in\mathcal{B}$;~~ b) It is the unique continuous trace on
$\mathcal{B}/\mathcal{B}^{-\infty}$.
\end{thm}

\section{The noncommutative residue $\widetilde{{\rm Wres}}[\pi^+(c(Z)\nabla_X^{S(TM)}\nabla_Y^{S(TM)}D^{-1})\circ\pi^+(D^{-2m})]$ on odd dimensional manifolds with boundary}
\label{section:3}
Firstly, we recall the definition of the Dirac operator. Let $M$ be an $n=(2m+1)$-dimensional oriented
compact spin Riemannian manifold with a Riemannian metric $g^{M}$ and let $\nabla^L$ be the Levi-Civita connection about $g^{M}$. In the fixed orthonormal frame $\{e_1,\cdots,e_n\}$, the connection matrix $(\omega_{s,t})$ is defined by
\begin{equation}
\label{a2}
\nabla^L(e_1,\cdots,e_n)= (e_1,\cdots,e_n)(\omega_{s,t}).
\end{equation}
\indent Let $c(e_i)$ denotes the Clifford action, which satisfies
\begin{align}
\label{a4}
&c(e_i)c(e_j)+c(e_j)c(e_i)=-2g^{M}(e_i,e_j).
\end{align}
By \cite{Y}, the Dirac operator is given by
\begin{align}
\label{a5}
&D=\sum^n_{i=1}c(e_i)[e_i-\frac{1}{4}\sum_{s,t}\omega_{s,t}
(e_i)c(e_s)c(e_t)].
\end{align}
\indent We define $\nabla_X^{S(TM)}:=X+\frac{1}{4}\sum_{ij}\langle\nabla_X^L{e_i},e_j\rangle c(e_i)c(e_j)$, which is a spin connection. Set\\ $A(X)=\frac{1}{4}\sum_{ij}\langle\nabla_X^L{e_i},e_j\rangle c(e_i)c(e_j)$, then
\begin{align}\label{ddd}
\nabla_X^{S(TM)}\nabla_Y^{S(TM)}&=[X+A(X)][Y+A(Y)]\nonumber\\
&=XY+X\cdot A(Y)+A(X)Y+A(X)A(Y)\nonumber\\
&=XY+X[A(Y)]+A(Y)X+A(X)Y+A(X)A(Y),
\end{align}
where $X=\sum_{j=1}^nX_j\partial_{x_j}, Y=\sum_{l=1}^nY_l\partial_{x_l}$ and $XY=\sum_{j,l=1}^n(X_jY_l\partial_{x_j}\partial_{x_l}+X_j\frac{\partial{Y_l}}{\partial_{x_j}}\partial_{x_l}).$ \\
\indent Let $g^{ij}=g(dx_{i},dx_{j})$, $\xi=\sum_{j}\xi_{j}dx_{j}$ and $\nabla^L_{\partial_{i}}\partial_{j}=\sum_{k}\Gamma_{ij}^{k}\partial_{k}$,  we denote that
\begin{align}
&\sigma_{i}=-\frac{1}{4}\sum_{s,t}\omega_{s,t}
(e_i)c(e_i)c(e_s)c(e_t)
;~~~\xi^{j}=g^{ij}\xi_{i};~~~~\Gamma^{k}=g^{ij}\Gamma_{ij}^{k};~~~~\sigma^{j}=g^{ij}\sigma_{i}.
\end{align}

And by $\sigma(\partial_{x_j})=\sqrt{-1}\xi_j$, we have the following lemmas.
\begin{lem}\label{lem3} The following identities hold:
\begin{align}
\label{b22}
\sigma_{0}(c(Z)\nabla_X^{S(TM)}\nabla_Y^{S(TM)})=&c(Z)X[A(Y)]+c(Z)A(X)A(Y);\nonumber\\
\sigma_{1}(c(Z)\nabla_X^{S(TM)}\nabla_Y^{S(TM)})=&\sqrt{-1}c(Z)\sum_{j,l=1}^nX_j\frac{\partial{Y_l}}{\partial_{x_j}}\xi_l+\sqrt{-1}c(Z)\sum_jA(Y)X_j\xi_j\nonumber\\
&+\sqrt{-1}c(Z)\sum_lA(X)Y_l\xi_l;\nonumber\\
\sigma_{2}(c(Z)\nabla_X^{S(TM)}\nabla_Y^{S(TM)})=&-c(Z)\sum_{j,l=1}^nX_jY_l\xi_j\xi_l.
\end{align}
\end{lem}
\begin{lem}\cite{Ka}\label{lem356} The following identities hold:
\begin{align}
\label{b22222}
&\sigma_{-2}(D^{-2})=|\xi|^{-2};\nonumber\\
&\sigma_{-3}(D^{-2})=-\sqrt{-1}|\xi|^{-4}\xi_k(\Gamma^k-2\sigma^k)-\sqrt{-1}|\xi|^{-6}2\xi^j\xi_\alpha\xi_\beta\partial_jg^{\alpha\beta}.
\end{align}
\end{lem}
\begin{lem}\cite{Wa5}\label{lem9356} The following identities hold:
\begin{align}
\label{bii}
&\sigma_{-1}(D^{-1})=\frac{\sqrt{-1}c(\xi)}{|\xi|^{2}};\nonumber\\
&\sigma_{-2}(D^{-1})=\frac{c(\xi)\sigma_0(D)c(\xi)}{|\xi|^4}
+\frac{c(\xi)}{|\xi|^6}\sum_jc(\mathrm{d}x_j)(\partial_{x_j}(c(\xi))|\xi|^2-c(\xi)\partial_{x_j}(|\xi|^2)),
\end{align}
where $\sigma_0(D)=-\frac{1}{4}\sum_{s,t,i}\omega_{s,t}(e_i)c(e_i)c(e_s)c(e_t)$.
\end{lem}
\indent Write
 \begin{eqnarray}
D_x^{\alpha}&=(-i)^{|\alpha|}\partial_x^{\alpha};
~\sigma(D_t)=p_1+p_0;
~(\sigma(D_t)^{-1})=\sum^{\infty}_{j=1}q_{-j}.
\end{eqnarray}

\indent By the composition formula of pseudo-differential operators, we have
\begin{align}
1=\sigma(D\circ D^{-1})&=\sum_{\alpha}\frac{1}{\alpha!}\partial^{\alpha}_{\xi}[\sigma(D)]
D_x^{\alpha}[\sigma(D^{-1})]\nonumber\\
&=(p_1+p_0)(q_{-1}+q_{-2}+q_{-3}+\cdots)\nonumber\\
&~~~+\sum_j(\partial_{\xi_j}p_1+\partial_{\xi_j}p_0)(
D_{x_j}q_{-1}+D_{x_j}q_{-2}+D_{x_j}q_{-3}+\cdots)\nonumber\\
&=p_1q_{-1}+(p_1q_{-2}+p_0q_{-1}+\sum_j\partial_{\xi_j}p_1D_{x_j}q_{-1})+\cdots,
\end{align}
so
\begin{equation}
q_{-1}=p_1^{-1};~q_{-2}=-p_1^{-1}[p_0p_1^{-1}+\sum_j\partial_{\xi_j}p_1D_{x_j}(p_1^{-1})].
\end{equation}

Then, we have the following lemma.
\begin{lem}\label{lema}The following identities hold:
\begin{align}\label{mki}
\sigma_{0}(c(Z)\nabla_X^{S(TM)}\nabla_Y^{S(TM)}D^{-2})&=-\sum_{j,l=1}^nX_jY_l\xi_j\xi_lc(Z)|\xi|^{-2};\nonumber\\
\sigma_{1}(c(Z)\nabla_X^{S(TM)}\nabla_Y^{S(TM)}D^{-1})&=-i\sum_{j,l=1}^nX_jY_l\xi_j\xi_lc(Z)c(\xi)|\xi|^{-2}.
\end{align}
\end{lem}

In this section, we compute the noncommutative residue $\widetilde{{\rm Wres}}[\pi^+(c(Z)\nabla_X^{S(TM)}\nabla_Y^{S(TM)}D^{-1})\circ\pi^+(D^{-2m})]$ on $(2m+1)$-dimensional oriented compact spin manifolds with boundary and get a general Dabrowski-Sitarz-Zalecki type theorem in this case.

Similar to \cite{Wa3}, by $\int_{S(\xi)}\xi_{j_{1}}\cdots\xi_{j_{l}}\texttt{d}(\xi)=0$ for odd $l$, we can compute the noncommutative residue
\begin{align}
\label{b141}
\widetilde{{\rm Wres}}[\pi^+(c(Z)\nabla_X^{S(TM)}\nabla_Y^{S(TM)}D^{-1})\circ\pi^+(D^{-2m})]=\int_{\partial M}\Phi,
\end{align}
where
\begin{align}
\label{qqq15}
\Phi =&\int_{|\xi'|=1}\int^{+\infty}_{-\infty}\sum^{\infty}_{j, k=0}\sum\frac{(-i)^{|\alpha|+j+k+1}}{\alpha!(j+k+1)!}
\times {\rm trace}_{S(TM)}[\partial^j_{x_n}\partial^\alpha_{\xi'}\partial^k_{\xi_n}\sigma^+_{r}(c(Z)\nabla_X^{S(TM)}\nabla_Y^{S(TM)}D^{-1})(x',0,\xi',\xi_n)\nonumber\\
&\times\partial^\alpha_{x'}\partial^{j+1}_{\xi_n}\partial^k_{x_n}\sigma_{l}(D^{-2m})(x',0,\xi',\xi_n)]d\xi_n\sigma(\xi')dx',
\end{align}
and the sum is taken over $r+l-k-j-|\alpha|-1=-(2m+1),~~r\leq 1,~~l\leq -2m$.

When $n=2m+1$ is odd, then ${\rm tr}_{S(TM)}[{\rm \texttt{id}}]={\rm dim}(\wedge^*(\mathbb{R}^2))=2^{\frac{n-1}{2}}=2^m$, where ${\rm tr}$ as shorthand of ${\rm trace}$, the sum is taken over $
r+l-k-j-|\alpha|=-2m,~~r\leq 1,~~l\leq -2m,$ then we have the following five cases:\\

\noindent  {\bf Case (a-I)}~$r=1,~l=-2m,~k=j=0,~|\alpha|=1$.\\

\noindent By (\ref{qqq15}), we get
\begin{equation}
\label{b24}
\Phi_1=-\int_{|\xi'|=1}\int^{+\infty}_{-\infty}\sum_{|\alpha|=1}
 {\rm tr}[\partial^\alpha_{\xi'}\pi^+_{\xi_n}\sigma_{1}(c(Z)\nabla_X^{S(TM)}\nabla_Y^{S(TM)}D^{-1})\times
 \partial^\alpha_{x'}\partial_{\xi_n}\sigma_{-2m}(D^{-2m})](x_0)d\xi_n\sigma(\xi')dx'.
\end{equation}
By Lemma 2.2 in \cite{Wa3}, for $i<n$,
\begin{equation}
\label{b25}
\partial_{x_i}\sigma_{-2m}(D^{-2m})(x_0)=\partial_{x_i}(|\xi|^2)^{-m}(x_0)=-m(|\xi|^2)^{-m-1}\partial_{x_i}(|\xi|^2)(x_0)=0,
\end{equation}
\noindent so $\Phi_1=0$.\\

 \noindent  {\bf Case (a-II)}~$r=1,~l=-2m,~k=|\alpha|=0,~j=1$.\\

\noindent By (\ref{qqq15}), we get
\begin{align}
\label{b26}
\Phi_2=&-\frac{1}{2}\int_{|\xi'|=1}\int^{+\infty}_{-\infty} {\rm
tr} [\partial_{x_n}\pi^+_{\xi_n}\sigma_{1}(c(Z)\nabla_X^{S(TM)}\nabla_Y^{S(TM)}D^{-1})\times
\partial_{\xi_n}^2\sigma_{-2m}(D^{-2m})](x_0)d\xi_n\sigma(\xi')dx'\nonumber\\
=&-\frac{1}{2}\int_{|\xi'|=1}\int^{+\infty}_{-\infty} {\rm
tr} [\partial_{\xi_n}^2\partial_{x_n}\pi^+_{\xi_n}\sigma_{1}(c(Z)\nabla_X^{S(TM)}\nabla_Y^{S(TM)}D^{-1})\times
\sigma_{-2m}(D^{-2m})](x_0)d\xi_n\sigma(\xi')dx'.
\end{align}
\noindent By (\ref{b22222}), we have\\
\begin{equation}\label{b237}
\sigma_{-2m}(D^{-2m})(x_0)=|\xi|^{-2m}=\frac{1}{(1+\xi_n^2)^m}.
\end{equation}
By Lemma \ref{lema}, we have\\
\begin{align}\label{b27}
\partial_{x_n}\sigma_{1}(c(Z)\nabla_X^{S(TM)}\nabla_Y^{S(TM)}D^{-1})(x_0)=&\partial_{x_n}\bigg(-i\sum_{j,l=1}^nX_jY_l\xi_j\xi_lc(Z)c(\xi)|\xi|^{-2}\bigg)(x_0)\nonumber\\
=&-i\sum_{j,l=1}^n\frac{\partial{X_j}}{\partial {x_n}}Y_l\xi_j\xi_l\frac{c(Z)c(\xi)}{|\xi|^2}-i\sum_{j,l=1}^nX_j\frac{\partial{Y_l}}{\partial {x_n}}\xi_j\xi_l\frac{c(Z)c(\xi)}{|\xi|^2}\nonumber\\
&-i\sum_{j,l=1}^nX_jY_l\xi_j\xi_l\frac{\partial_{x_n}(c(Z))c(\xi)}{|\xi|^2}-i\sum_{j,l=1}^nX_jY_l\xi_j\xi_l\frac{c(Z)\partial_{x_n}[c(\xi')](x_0)}{|\xi|^2}\nonumber\\
&+i\sum_{j,l=1}^nX_jY_l\xi_j\xi_l\frac{c(Z)c(\xi)|\xi'|^2h'(0)}{|\xi|^4}.
\end{align}
If we omit some items that have no contribution for computing ${\bf \Phi_2}$. Then, we have
\begin{align}\label{b28}
&\partial_{x_n}\pi^+_{\xi_n}\sigma_{1}(c(Z)\nabla_X^{S(TM)}\nabla_Y^{S(TM)}D^{-1})(x_0)\nonumber\\
=&\pi^+_{\xi_n}\partial_{x_n}\sigma_{1}(c(Z)\nabla_X^{S(TM)}\nabla_Y^{S(TM)}D^{-1})(x_0)\nonumber\\
=&-\sum_{j,l=1}^{n-1}\bigg(\frac{\partial{X_j}}{\partial {x_n}}Y_l+X_j\frac{\partial{Y_l}}{\partial {x_n}}\bigg)\xi_j\xi_l\bigg[\frac{c(Z)c(\xi')}{2(\xi_n-i)}+\frac{ic(Z)c(dx_n)}{2(\xi_n-i)}\bigg]\nonumber\\
&-i\bigg(\frac{\partial{X_n}}{\partial {x_n}}Y_n+X_n\frac{\partial{Y_n}}{\partial {x_n}}\bigg)\bigg[\frac{ic(Z)c(\xi')}{2(\xi_n-i)}-\frac{c(Z)c(dx_n)}{2(\xi_n-i)}\bigg]\nonumber\\
&-\sum_{j,l=1}^{n-1}X_jY_l\xi_j\xi_l\bigg\{\frac{\partial_{x_n}(c(Z))c(\xi')}{2(\xi_n-i)}+\frac{i\partial_{x_n}(c(Z))c(dx_n)}{2(\xi_n-i)}+\frac{c(Z)\partial_{x_n}[c(\xi')](x_0)}{2(\xi_n-i)}\nonumber\\
&+ih'(0)|\xi'|^2\bigg[\frac{(i\xi_n+2)c(Z)c(\xi')}{4(\xi_n-i)^2}+\frac{ic(Z)c(dx_n)}{4(\xi_n-i)^2}\bigg]\bigg\}\nonumber\\
&-iX_nY_n\bigg\{\frac{i\partial_{x_n}(c(Z))c(\xi')}{2(\xi_n-i)}-\frac{\partial_{x_n}(c(Z))c(dx_n)}{2(\xi_n-i)}+\frac{ic(Z)\partial_{x_n}[c(\xi')](x_0)}{2(\xi_n-i)}\nonumber\\
&-h'(0)|\xi'|^2\bigg[\frac{-i\xi_nc(Z)c(\xi')}{4(\xi_n-i)^2}+\frac{(2\xi_n-i)c(Z)c(dx_n)}{4(\xi_n-i)^2}\bigg]\bigg\}.
\end{align}
By further calculation, we have
\begin{align}\label{b28}
&\partial_{\xi_n}^2\partial_{x_n}\pi^+_{\xi_n}\sigma_{1}(c(Z)\nabla_X^{S(TM)}\nabla_Y^{S(TM)}D^{-1})(x_0)\nonumber\\
=&-\sum_{j,l=1}^{n-1}\bigg(\frac{\partial{X_j}}{\partial {x_n}}Y_l+X_j\frac{\partial{Y_l}}{\partial {x_n}}\bigg)\xi_j\xi_l\bigg[\frac{c(Z)c(\xi')}{(\xi_n-i)^3}+\frac{ic(Z)c(dx_n)}{(\xi_n-i)^3}\bigg]\nonumber\\
&-i\bigg(\frac{\partial{X_n}}{\partial {x_n}}Y_n+X_n\frac{\partial{Y_n}}{\partial {x_n}}\bigg)\bigg[\frac{ic(Z)c(\xi')}{(\xi_n-i)^3}-\frac{c(Z)c(dx_n)}{(\xi_n-i)^3}\bigg]\nonumber\\
&-\sum_{j,l=1}^{n-1}X_jY_l\xi_j\xi_l\bigg\{\frac{\partial_{x_n}(c(Z))c(\xi')}{(\xi_n-i)^3}+\frac{i\partial_{x_n}(c(Z))c(dx_n)}{(\xi_n-i)^3}+\frac{c(Z)\partial_{x_n}[c(\xi')](x_0)}{(\xi_n-i)^3}\nonumber\\
&+ih'(0)|\xi'|^2\bigg[\frac{(i\xi_n+4)c(Z)c(\xi')}{2(\xi_n-i)^4}+\frac{3ic(Z)c(dx_n)}{2(\xi_n-i)^4}\bigg]\bigg\}\nonumber\\
&-iX_nY_n\bigg\{\frac{i\partial_{x_n}(c(Z))c(\xi')}{(\xi_n-i)^3}-\frac{\partial_{x_n}(c(Z))c(dx_n)}{(\xi_n-i)^3}+\frac{ic(Z)\partial_{x_n}[c(\xi')](x_0)}{(\xi_n-i)^3}\nonumber\\
&-h'(0)|\xi'|^2\bigg[\frac{(-i\xi_n+2)c(Z)c(\xi')}{2(\xi_n-i)^4}+\frac{(2\xi_n+i)c(Z)c(dx_n)}{2(\xi_n-i)^4}\bigg]\bigg\}.
\end{align}
By the relation of the Clifford action and ${\rm tr}{(AB)}={\rm tr }{(BA)}$,  we have the equalities:
\begin{align}\label{a49}
&{\rm tr}[c(Z)
c(\xi')]=-g(Z,\xi'){\rm tr}[\texttt{id}];~~
{\rm tr}[c(Z)c(dx_n)]=-Z_n{\rm tr}[\texttt{id}];~~{\rm tr}[c(Z)
\partial_{x_n}[c(\xi')]]=-\frac{1}{2}h'(0)g(Z,\xi'){\rm tr}[\texttt{id}];\nonumber\\
&{\rm tr }[\partial_{x_n}(c(Z))c(\xi')]=-\partial_{x_n}(g(Z,\xi')){\rm tr}[\texttt{id}]+\frac{1}{2}h'(0)g(Z,\xi'){\rm tr}[\texttt{id}];
~~{\rm tr}[\partial_{x_n}(c(Z))c(dx_n)]=-\partial_{x_n}(Z_n){\rm tr}[\texttt{id}].
\end{align}
Let $X=X^T+X_n\partial_n,~Y=Y^T+Y_n\partial_n,$ then we have $\sum_{j=1}^{n-1}X_jY_j(x_0)=g(X^T,Y^T)(x_0)$. And by $\int_{S^{n-2}}\xi_j\xi_l\sigma(\xi')=\frac{1}{n-1}\delta_{jl}Vol(S^{n-2})$,  where $Vol(S^{n-2})$ is the canonical volume of $S^{n-2}$.\\
We note that $i<n,~\int_{|\xi'|=1}\{\xi_{i_{1}}\xi_{i_{2}}\cdots\xi_{i_{2q+1}}\}\sigma(\xi')=0$,
so $g(Z,\xi')$ and $\partial_{x_n}(g(Z,\xi'))$ have no contribution for computing ${\bf \Phi_2}$. Then, we have
\begin{align}\label{35}
&{\rm tr} [\partial_{\xi_n}^2\partial_{x_n}\pi^+_{\xi_n}\sigma_{1}(c(Z)\nabla_X^{S(TM)}\nabla_Y^{S(TM)}D^{-1})\times
\sigma_{-2m}(D^{-2m})](x_0)|_{|\xi'|=1}\nonumber\\
=&\sum_{j,l=1}^{n-1}\bigg(\frac{\partial{X_j}}{\partial {x_n}}Y_l+X_j\frac{\partial{Y_l}}{\partial {x_n}}\bigg)\xi_j\xi_l\frac{i}{(\xi_n-i)^{m+3}(\xi_n+i)^{m}}Z_n{\rm tr}[\texttt{id}]\nonumber\\
&-\bigg(\frac{\partial{X_n}}{\partial {x_n}}Y_n+X_n\frac{\partial{Y_n}}{\partial {x_n}}\bigg)\frac{i}{(\xi_n-i)^{m+3}(\xi_n+i)^{m}}Z_n{\rm tr}[\texttt{id}]\nonumber\\
&+\sum_{j,l=1}^{n-1}X_jY_l\xi_j\xi_l\frac{i}{(\xi_n-i)^{m+3}(\xi_n+i)^{m}}\partial_{x_n}(Z_n){\rm tr}[\texttt{id}]\nonumber\\
&-\sum_{j,l=1}^{n-1}X_jY_l\xi_j\xi_l\frac{3}{2(\xi_n-i)^{m+4}(\xi_n+i)^{m}}Z_nh'(0){\rm tr}[\texttt{id}]\nonumber\\
&-X_nY_n\frac{i}{(\xi_n-i)^{m+3}(\xi_n+i)^{m}}\partial_{x_n}(Z_n){\rm tr}[\texttt{id}]\nonumber\\
&-X_nY_n\frac{2i\xi_n-1}{2(\xi_n-i)^{m+4}(\xi_n+i)^{m}}Z_nh'(0){\rm tr}[\texttt{id}].
\end{align}
Therefore, we get
\begin{align}\label{35kkk}
\Phi_2
=&-\frac{1}{2}\int_{|\xi'|=1}\int^{+\infty}_{-\infty} \sum_{j,l=1}^{n-1}\bigg(\frac{\partial{X_j}}{\partial {x_n}}Y_l+X_j\frac{\partial{Y_l}}{\partial {x_n}}\bigg)\xi_j\xi_l\frac{i}{(\xi_n-i)^{m+3}(\xi_n+i)^{m}}Z_n{\rm tr}[\texttt{id}]d\xi_n\sigma(\xi')dx'\nonumber\\
&+\frac{1}{2}\int_{|\xi'|=1}\int^{+\infty}_{-\infty} \bigg(\frac{\partial{X_n}}{\partial {x_n}}Y_n+X_n\frac{\partial{Y_n}}{\partial {x_n}}\bigg)\frac{i}{(\xi_n-i)^{m+3}(\xi_n+i)^{m}}Z_n{\rm tr}[\texttt{id}]d\xi_n\sigma(\xi')dx'\nonumber\\
&-\frac{1}{2}\int_{|\xi'|=1}\int^{+\infty}_{-\infty} \sum_{j,l=1}^{n-1}X_jY_l\xi_j\xi_l\frac{i}{(\xi_n-i)^{m+3}(\xi_n+i)^{m}}\partial_{x_n}(Z_n){\rm tr}[\texttt{id}]d\xi_n\sigma(\xi')dx'\nonumber\\
&+\frac{1}{2}\int_{|\xi'|=1}\int^{+\infty}_{-\infty} \sum_{j,l=1}^{n-1}X_jY_l\xi_j\xi_l\frac{3}{2(\xi_n-i)^{m+4}(\xi_n+i)^{m}}Z_nh'(0){\rm tr}[\texttt{id}]d\xi_n\sigma(\xi')dx'\nonumber\\
&+\frac{1}{2}\int_{|\xi'|=1}\int^{+\infty}_{-\infty} X_nY_n\frac{i}{(\xi_n-i)^{m+3}(\xi_n+i)^{m}}\partial_{x_n}(Z_n){\rm tr}[\texttt{id}]d\xi_n\sigma(\xi')dx'\nonumber\\
&+\frac{1}{2}\int_{|\xi'|=1}\int^{+\infty}_{-\infty} X_nY_n\frac{2i\xi_n-1}{2(\xi_n-i)^{m+4}(\xi_n+i)^{m}}Z_nh'(0){\rm tr}[\texttt{id}]d\xi_n\sigma(\xi')dx'\nonumber\\
=&-\frac{i}{4m}Vol(S^{n-2})\partial_{x_n}\Big(g(X^T,Y^T)Z_n\Big)2^m\int_{\Gamma^{+}}\frac{1}{(\xi_n-i)^{m+3}(\xi_n+i)^{m}}d\xi_ndx'\nonumber\\
&+\frac{i}{2}Vol(S^{n-2})\partial_{x_n}\Big((X_nY_n)Z_n\Big)2^m\int_{\Gamma^{+}}\frac{1}{(\xi_n-i)^{m+3}(\xi_n+i)^{m}}d\xi_ndx'\nonumber\\
&+\frac{1}{8m}Vol(S^{n-2})g(X^T,Y^T)Z_nh'(0)2^m\int_{\Gamma^{+}}\frac{3}{(\xi_n-i)^{m+4}(\xi_n+i)^{m}}d\xi_ndx'\nonumber\\
&+\frac{1}{4}Vol(S^{n-2})X_nY_nZ_nh'(0)2^m\int_{\Gamma^{+}}\frac{2i\xi_n-1}{(\xi_n-i)^{m+4}(\xi_n+i)^{m}}d\xi_ndx'\nonumber\\
=&-\frac{i}{4}Vol(S^{n-2})\partial_{x_n}\Big(g(X^T,Y^T)Z_n\Big)2^m\frac{2\pi i}{(m+2)!}\left[\frac{1}{m(\xi_n+i)^{m}}\right]^{(m+2)}\bigg|_{\xi_n=i}dx'\nonumber\\
&+\frac{i}{2}Vol(S^{n-2})\partial_{x_n}\Big((X_nY_n)Z_n\Big)2^m\frac{2\pi i}{(m+2)!}\left[\frac{1}{(\xi_n+i)^{m}}\right]^{(m+2)}\bigg|_{\xi_n=i}dx'\nonumber\\
&+\frac{1}{8}Vol(S^{n-2})g(X^T,Y^T)Z_nh'(0)2^m\frac{2\pi i}{(m+3)!}\left[\frac{3}{m(\xi_n+i)^{m}}\right]^{(m+3)}\bigg|_{\xi_n=i}dx'\nonumber\\
&+\frac{1}{4}Vol(S^{n-2})X_nY_nZ_nh'(0)2^m\frac{2\pi i}{(m+3)!}\left[\frac{2i\xi_n-1}{(\xi_n+i)^{m}}\right]^{(m+3)}\bigg|_{\xi_n=i}dx'\nonumber\\
:=&\bigg\{\partial_{x_n}\Big(g(X^T,Y^T)Z_n\Big)\frac{\pi}{2(m+2)!}A_0-\partial_{x_n}\Big((X_nY_n)Z_n\Big)\frac{\pi}{(m+2)!}A_1\nonumber\\
&+Z_nh'(0)\Big(g(X^T,Y^T)\frac{\pi i}{4(m+3)!}A_2+X_nY_n\frac{\pi i}{2(m+3)!}A_3\Big)\bigg\}Vol(S^{n-2})2^mdx',
\end{align}
where
\begin{align}
A_0&=\left[\frac{1}{m(\xi_n+i)^{m}}\right]^{m+2}\bigg|_{\xi_n=i}=-\frac{(m+1)(m+2)\cdot\cdot\cdot(2m+1)}{2^{2m+2}};\nonumber\\
A_1&=\left[\frac{1}{(\xi_n+i)^{m}}\right]^{(m+2)}\bigg|_{\xi_n=i}=-\frac{m(m+1)\cdot\cdot\cdot(2m+1)}{2^{2m+2}};\nonumber\\
A_2&=\left[\frac{3}{m(\xi_n+i)^{m}}\right]^{(m+3)}\bigg|_{\xi_n=i}=\frac{3(m+1)(m+2)\cdot\cdot\cdot(2m+2)}{2^{2m+3}i};\nonumber\\
A_3&=\left[\frac{2i\xi_n-1}{(\xi_n+i)^{m}}\right]^{(m+3)}\bigg|_{\xi_n=i}=\frac{(3-m)m(m+1)\cdot\cdot\cdot(2m+1)}{2^{2m+2}i}.
\end{align}\\

\noindent  {\bf Case (a-III)}~$r=1,~l=-2m,~j=|\alpha|=0,~k=1$.\\

\noindent By (\ref{qqq15}), we get
\begin{align}\label{36}
\Phi_3&=-\frac{1}{2}\int_{|\xi'|=1}\int^{+\infty}_{-\infty}
{\rm tr} [\partial_{\xi_n}\pi^+_{\xi_n}\sigma_{1}(c(Z)\nabla_X^{S(TM)}\nabla_Y^{S(TM)}D^{-1})\times
\partial_{\xi_n}\partial_{x_n}\sigma_{-2m}(D^{-2m})](x_0)d\xi_n\sigma(\xi')dx'\nonumber\\
&=\frac{1}{2}\int_{|\xi'|=1}\int^{+\infty}_{-\infty}
{\rm tr} [\partial_{\xi_n}^2\pi^+_{\xi_n}\sigma_{1}(c(Z)\nabla_X^{S(TM)}\nabla_Y^{S(TM)}D^{-1})\times
\partial_{x_n}\sigma_{-2m}(D^{-2m})](x_0)d\xi_n\sigma(\xi')dx'.
\end{align}
By (3.21) in \cite{Wa6}, we have
\begin{eqnarray}\label{37}
\partial_{x_n} \sigma_{-2m}(D^{-2m})(x_0)
=\partial_{x_n}(|\xi|^2)^{-m} (x_0)
=\frac{-mh'(0)|\xi'|^2}{(1+\xi_n^2)^{m+1}}.
\end{eqnarray}
By Lemma \ref{lema}, we have
\begin{align}\label{aa388}
\pi^+_{\xi_n}\sigma_{1}(c(Z)\nabla_X^{S(TM)}\nabla_Y^{S(TM)}D^{-1})(x_0)=&-\sum_{j,l=1}^{n-1}X_jY_l\xi_j\xi_l\bigg[\frac{c(Z)c(\xi')}{2(\xi_n-i)}+\frac{ic(Z)c(dx_n)}{2(\xi_n-i)}\bigg]\nonumber\\
&+X_nY_n\bigg[\frac{c(Z)c(\xi')}{2(\xi_n-i)}+\frac{ic(Z)c(dx_n)}{2(\xi_n-i)}\bigg]\nonumber\\
&-\sum_{j=1}^{n-1}(X_jY_n+X_nY_j)\xi_j\bigg[\frac{ic(Z)c(\xi')}{2(\xi_n-i)}-\frac{c(Z)c(dx_n)}{2(\xi_n-i)}\bigg].
\end{align}
Then, we get
\begin{align}\label{mmmmm}
\partial_{\xi_n}^2\pi^+_{\xi_n}\sigma_{1}(c(Z)\nabla_X^{S(TM)}\nabla_Y^{S(TM)}D^{-1})(x_0)=&-\sum_{j,l=1}^{n-1}X_jY_l\xi_j\xi_l\bigg[\frac{c(Z)c(\xi')}{(\xi_n-i)^3}+\frac{ic(Z)c(dx_n)}{(\xi_n-i)^3}\bigg]\nonumber\\
&+X_nY_n\bigg[\frac{c(Z)c(\xi')}{(\xi_n-i)^3}+\frac{ic(Z)c(dx_n)}{(\xi_n-i)^3}\bigg]\nonumber\\
&-\sum_{j=1}^{n-1}(X_jY_n+X_nY_j)\xi_j\bigg[\frac{ic(Z)c(\xi')}{(\xi_n-i)^3}-\frac{c(Z)c(dx_n)}{(\xi_n-i)^3}\bigg].
\end{align}
We omit some items that have no contribution for computing ${\bf \Phi_3}$. Then, we have
\begin{align}\label{35}
&{\rm tr} [\partial_{\xi_n}^2\pi^+_{\xi_n}\sigma_{1}(c(Z)\nabla_X^{S(TM)}\nabla_Y^{S(TM)}D^{-1})\times
\partial_{x_n}\sigma_{-2m}(D^{-2m})](x_0)|_{|\xi'|=1}\nonumber\\
=&-\sum_{j,l=1}^{n-1}X_jY_l\xi_j\xi_l\frac{mi}{(\xi_n-i)^{m+4}(\xi_n+i)^{m+1}}Z_nh'(0){\rm tr}[\texttt{id}]\nonumber\\
&+X_nY_n\frac{mi}{(\xi_n-i)^{m+4}(\xi_n+i)^{m+1}}Z_nh'(0){\rm tr}[\texttt{id}].
\end{align}
Therefore, we get
\begin{align}\label{41}
\Phi_3=&-\frac{1}{2}\int_{|\xi'|=1}\int^{+\infty}_{-\infty} \sum_{j,l=1}^{n-1}X_jY_l\xi_j\xi_l\frac{mi}{(\xi_n-i)^{m+4}(\xi_n+i)^{m+1}}Z_nh'(0){\rm tr}[\texttt{id}]d\xi_n\sigma(\xi')dx'\nonumber\\
&+\frac{1}{2}\int_{|\xi'|=1}\int^{+\infty}_{-\infty} X_nY_n\frac{mi}{(\xi_n-i)^{m+4}(\xi_n+i)^{m+1}}Z_nh'(0){\rm tr}[\texttt{id}]d\xi_n\sigma(\xi')dx'\nonumber\\
=&-\frac{i}{4}Vol(S^{n-2})g(X^T,Y^T)Z_nh'(0)2^m\int_{\Gamma^{+}} \frac{1}{(\xi_n-i)^{m+4}(\xi_n+i)^{m+1}}d\xi_ndx'\nonumber\\
&+\frac{i}{2}Vol(S^{n-2})X_nY_nZ_nh'(0)2^m\int_{\Gamma^{+}} \frac{m}{(\xi_n-i)^{m+4}(\xi_n+i)^{m+1}}d\xi_ndx'\nonumber\\
=&-\frac{i}{4}Vol(S^{n-2})g(X^T,Y^T)Z_nh'(0)2^m\frac{2\pi i}{(m+3)!}\left[\frac{1}{(\xi_n+i)^{m+1}}\right]^{(m+3)}\bigg|_{\xi_n=i}dx'\nonumber\\
&+\frac{i}{2}Vol(S^{n-2})X_nY_nZ_nh'(0)2^m\frac{2\pi i}{(m+3)!}\left[\frac{m}{(\xi_n+i)^{m+1}}\right]^{(m+3)}\bigg|_{\xi_n=i}dx'\nonumber\\
:=&\Big(g(X^T,Y^T)\frac{\pi}{2(m+3)!}B_0-X_nY_n\frac{\pi}{(m+3)!}B_1\Big)Vol(S^{n-2})Z_nh'(0)2^mdx',
\end{align}
where
\begin{align}
B_0&=\left[\frac{1}{(\xi_n+i)^{m+1}}\right]^{(m+3)}\bigg|_{\xi_n=i}=-\frac{(m+1)(m+2)\cdot\cdot\cdot(2m+3)}{2^{2m+4}};\nonumber\\
B_1&=\left[\frac{m}{(\xi_n+i)^{m+1}}\right]^{(m+3)}\bigg|_{\xi_n=i}=-\frac{m(m+1)\cdot\cdot\cdot(2m+3)}{2^{2m+4}}.
\end{align}\\

\noindent  {\bf Case (a-IV)}~$r=1,~l=-(2m+1),~k=j=|\alpha|=0$.\\

\noindent By (\ref{qqq15}), we get
\begin{align}\label{42}
\Phi_4&=-i\int_{|\xi'|=1}\int^{+\infty}_{-\infty}{\rm tr} [\pi^+_{\xi_n}\sigma_{1}(c(Z)\nabla_X^{S(TM)}\nabla_Y^{S(TM)}D^{-1})\times
\partial_{\xi_n}\sigma_{-(2m+1)}(D^{-2m})](x_0)d\xi_n\sigma(\xi')dx'\nonumber\\
&=i\int_{|\xi'|=1}\int^{+\infty}_{-\infty}{\rm tr} [\partial_{\xi_n}\pi^+_{\xi_n}\sigma_{1}(c(Z)\nabla_X^{S(TM)}\nabla_Y^{S(TM)}D^{-1})\times
\sigma_{-(2m+1)}(D^{-2m})](x_0)d\xi_n\sigma(\xi')dx'.
\end{align}
By (3.30) in \cite{Wa6}, we have
\begin{align}\label{43}
\sigma_{-(2m+1)}(D^{-2m})(x_0)
=&m(1+\xi_n^2)^{-m+1}\Big[\frac{-i}{(1+\xi_n^2)^2}\times\frac{2m+3}{2}h'(0)\xi_n-\frac{2ih'(0)\xi_n}{(1+\xi_n^2)^3} \Big]\nonumber\\
&+ih'(0)(-m^2+m)\xi_n(1+\xi_n^2)^{-m-2}\nonumber\\
=&-\frac{(2m^2+3m)ih'(0)\xi_n}{2(1+\xi_n^2)^{m+1}}-\frac{(m^2+m)ih'(0)\xi_n}{(1+\xi_n^2)^{m+2}}.
\end{align}
By (\ref{aa388}), we have
\begin{align}\label{45}
\partial_{\xi_n}\pi^+_{\xi_n}\sigma_{1}(c(Z)\nabla_X^{S(TM)}\nabla_Y^{S(TM)}D^{-1})(x_0)=&\sum_{j,l=1}^{n-1}X_jY_l\xi_j\xi_l\bigg[\frac{c(Z)c(\xi')}{2(\xi_n-i)^2}+\frac{ic(Z)c(dx_n)}{2(\xi_n-i)^2}\bigg]\nonumber\\
&-X_nY_n\bigg[\frac{c(Z)c(\xi')}{2(\xi_n-i)^2}+\frac{ic(Z)c(dx_n)}{2(\xi_n-i)^2}\bigg]\nonumber\\
&+\sum_{j=1}^{n-1}(X_jY_n+X_nY_j)\xi_j\bigg[\frac{ic(Z)c(\xi')}{2(\xi_n-i)^2}-\frac{c(Z)c(dx_n)}{2(\xi_n-i)^2}\bigg].
\end{align}
We omit some items that have no contribution for computing ${\bf \Phi_4}$.
Then, we have
\begin{align}\label{35}
&{\rm tr} [\partial_{\xi_n}\pi^+_{\xi_n}\sigma_{1}(c(Z)\nabla_X^{S(TM)}\nabla_Y^{S(TM)}D^{-1})\times
\sigma_{-(2m+1)}(D^{-2m})](x_0)|_{|\xi'|=1}\nonumber\\
=&-\sum_{j,l=1}^{n-1}X_jY_l\xi_j\xi_l\bigg[\frac{(2m^2+3m)\xi_n}{4(\xi_n-i)^{m+3}(\xi_n+i)^{m+1}}+\frac{(m^2+m)\xi_n}{2(\xi_n-i)^{m+4}(\xi_n+i)^{m+2}}\bigg]Z_nh'(0){\rm tr}[\texttt{id}]\nonumber\\
&+X_nY_n\bigg[\frac{(2m^2+3m)\xi_n}{4(\xi_n-i)^{m+3}(\xi_n+i)^{m+1}}+\frac{(m^2+m)\xi_n}{2(\xi_n-i)^{m+4}(\xi_n+i)^{m+2}}\bigg]Z_nh'(0){\rm tr}[\texttt{id}].
\end{align}
Therefore, we get
\begin{align}\label{39}
\Phi_4=&-i\int_{|\xi'|=1}\int^{+\infty}_{-\infty}\sum_{j,l=1}^{n-1}X_jY_l\xi_j\xi_l\frac{(2m^2+3m)\xi_n}{4(\xi_n-i)^{m+3}(\xi_n+i)^{m+1}}Z_nh'(0){\rm tr}[\texttt{id}]d\xi_n\sigma(\xi')dx'\nonumber\\
&-i\int_{|\xi'|=1}\int^{+\infty}_{-\infty}\sum_{j,l=1}^{n-1}X_jY_l\xi_j\xi_l\frac{(m^2+m)\xi_n}{2(\xi_n-i)^{m+4}(\xi_n+i)^{m+2}}Z_nh'(0){\rm tr}[\texttt{id}]d\xi_n\sigma(\xi')dx'\nonumber\\
&+i\int_{|\xi'|=1}\int^{+\infty}_{-\infty}X_nY_n\frac{(2m^2+3m)\xi_n}{4(\xi_n-i)^{m+3}(\xi_n+i)^{m+1}}Z_nh'(0){\rm tr}[\texttt{id}]d\xi_n\sigma(\xi')dx'\nonumber\\
&+i\int_{|\xi'|=1}\int^{+\infty}_{-\infty}X_nY_n\frac{(m^2+m)\xi_n}{2(\xi_n-i)^{m+4}(\xi_n+i)^{m+2}}Z_nh'(0){\rm tr}[\texttt{id}]d\xi_n\sigma(\xi')dx'\nonumber\\
=&-\frac{i}{8}Vol(S^{n-2})g(X^T,Y^T)Z_nh'(0)2^m\int_{\Gamma^{+}} \frac{(2m+3)\xi_n}{(\xi_n-i)^{m+3}(\xi_n+i)^{m+1}}d\xi_ndx'\nonumber\\
&-\frac{i}{4}Vol(S^{n-2})g(X^T,Y^T)Z_nh'(0)2^m\int_{\Gamma^{+}} \frac{(m+1)\xi_n}{(\xi_n-i)^{m+4}(\xi_n+i)^{m+2}}d\xi_ndx'\nonumber\\
&+\frac{i}{4}Vol(S^{n-2})X_nY_nZ_nh'(0)2^m\int_{\Gamma^{+}} \frac{(2m^2+3m)\xi_n}{(\xi_n-i)^{m+3}(\xi_n+i)^{m+1}}d\xi_ndx'\nonumber\\
&+\frac{i}{2}Vol(S^{n-2})X_nY_nZ_nh'(0)2^m\int_{\Gamma^{+}} \frac{(m^2+m)\xi_n}{(\xi_n-i)^{m+4}(\xi_n+i)^{m+2}}d\xi_ndx'\nonumber\\
=&-\frac{i}{8}Vol(S^{n-2})g(X^T,Y^T)Z_nh'(0)2^m\frac{2\pi i}{(m+2)!}\left[\frac{(2m+3)\xi_n}{(\xi_n+i)^{m+1}}\right]^{(m+2)}\bigg|_{\xi_n=i}dx'\nonumber\\
&-\frac{i}{4}Vol(S^{n-2})g(X^T,Y^T)Z_nh'(0)2^m\frac{2\pi i}{(m+3)!}\left[\frac{(m+1)\xi_n}{(\xi_n+i)^{m+2}}\right]^{(m+3)}\bigg|_{\xi_n=i}dx'\nonumber\\
&+\frac{i}{4}Vol(S^{n-2})X_nY_nZ_nh'(0)2^m\frac{2\pi i}{(m+2)!}\left[\frac{(2m^2+3m)\xi_n}{(\xi_n+i)^{m+1}}\right]^{(m+2)}\bigg|_{\xi_n=i}dx'\nonumber\\
&+\frac{i}{2}Vol(S^{n-2})X_nY_nZ_nh'(0)2^m\frac{2\pi i}{(m+3)!}\left[\frac{(m^2+m)\xi_n}{(\xi_n+i)^{m+2}}\right]^{(m+3)}\bigg|_{\xi_n=i}dx'\nonumber\\
:=&\bigg\{g(X^T,Y^T)\Big(\frac{\pi}{4(m+2)!}C_0+\frac{\pi}{2(m+3)!}C_1\Big)-X_nY_n\Big(\frac{\pi}{2(m+2)!}C_2\nonumber\\
&+\frac{\pi}{(m+3)!}C_3\Big)\bigg\}Vol(S^{n-2})Z_nh'(0)2^mdx',
\end{align}
where
\begin{align}
C_0&=\left[\frac{(2m+3)\xi_n}{(\xi_n+i)^{m+1}}\right]^{(m+2)}\bigg|_{\xi_n=i}=-\frac{(m+2)(m+3)\cdot\cdot\cdot(2m+3)}{2^{2m+3}};\nonumber\\
C_1&=\left[\frac{(m+1)\xi_n}{(\xi_n+i)^{m+2}}\right]^{(m+3)}\bigg|_{\xi_n=i}=\frac{(m+1)(m+2)\cdot\cdot\cdot(2m+3)}{2^{2m+4}};\nonumber\\
C_2&=\left[\frac{(2m^2+3m)\xi_n}{(\xi_n+i)^{m+1}}\right]^{(m+2)}\bigg|_{\xi_n=i}=-\frac{m(m+2)(m+3)\cdot\cdot\cdot(2m+3)}{2^{2m+3}};\nonumber\\
C_3&=\left[\frac{(m^2+m)\xi_n}{(\xi_n+i)^{m+2}}\right]^{(m+3)}\bigg|_{\xi_n=i}=\frac{m(m+1)\cdot\cdot\cdot(2m+3)}{2^{2m+4}}.
\end{align}\\

\noindent {\bf Case (a-V)}~$r=0,~l=-2m,~k=j=|\alpha|=0$.\\

\noindent By (\ref{qqq15}), we get
\begin{align}\label{61}
\Phi_5=-i\int_{|\xi'|=1}\int^{+\infty}_{-\infty}{\rm tr} [\pi^+_{\xi_n}\sigma_{0}(c(Z)\nabla_X^{S(TM)}\nabla_Y^{S(TM)}D^{-1})\times
\partial_{\xi_n}\sigma_{-2m}(D^{-2m})](x_0)d\xi_n\sigma(\xi')dx'.
\end{align}
By (3.33) in \cite{Wa6}, we have
\begin{equation}\label{62}
\partial_{\xi_n}\sigma_{-2m}(D^{-2m})(x_0)
=\partial_{\xi_n}(|\xi|^2)^{-m}(x_0)=\frac{-2m\xi_n}{(1+\xi_n^2)^{m+1}}.
\end{equation}
By Lemma \ref{lem3} and Lemma \ref{lem9356}, we have
\begin{align}\label{ggg}
\sigma_{0}(c(Z)\nabla_X^{S(TM)}\nabla_Y^{S(TM)}D^{-1})=&\sigma_{1}(c(Z)\nabla_X^{S(TM)}\nabla_Y^{S(TM)})\sigma_{-1}(D^{-1})+\sigma_{2}(c(Z)\nabla_X^{S(TM)}\nabla_Y^{S(TM)})\sigma_{-2}(D^{-1})\nonumber\\
&+\sum_{j=1}^n\partial_{\xi_j}\sigma_{2}(c(Z)\nabla_X^{S(TM)}\nabla_Y^{S(TM)})D_{x_j}[\sigma_{-1}(D^{-1})]\nonumber\\
:=&A^1+A^2+A^3,
\end{align}
where
\begin{align}
A^1(x_0)=&\sqrt{-1}c(Z)\bigg(\sum_{j,l=1}^nX_j\frac{\partial{Y_l}}{\partial_{x_j}}\xi_l+\sum_jA(Y)X_j\xi_j+\sum_lA(X)Y_l\xi_l\bigg)\frac{\sqrt{-1}c(\xi)}{|\xi|^{2}};\nonumber\\
A^2(x_0)=&-c(Z)\sum_{j,l=1}^{n}X_jY_l\xi_j\xi_l\bigg[\frac{c(\xi)\sigma_0(D)(x_0)c(\xi)}{|\xi|^4}
+\frac{c(\xi)}{|\xi|^6}\sum_jc(\mathrm{d}x_j)(\partial_{x_j}(c(\xi))|\xi|^2-c(\xi)\partial_{x_j}(|\xi|^2))\bigg];\nonumber\\
A^3(x_0)=&-c(Z)\bigg(\sum_{l=1}^{n-1}(X_nY_l+X_lY_n)\xi_l+2X_nY_n\xi_n\bigg)\bigg[\frac{\partial_{x_n}[c(\xi')](x_0)}{|\xi|^2}-\frac{h'(0)c(\xi)}{|\xi|^4}\bigg].
\end{align}
Firstly, the following results are obtained by further calculation of $A^1(x_0$)
\begin{align}
A^1(x_0)=&-\sum_{j=1}^n\sum_{l=1}^{n-1}X_j\frac{\partial {Y_l}}{\partial {x_j}}\xi_l\frac{c(Z)c(\xi)}{|\xi|^2}-\sum_{j=1}^nX_j\frac{\partial {Y_n}}{\partial {x_j}}\xi_n\frac{c(Z)c(\xi)}{|\xi|^2}-c(Z)\sum_{j=1}^{n-1}A(Y)X_j\xi_j\frac{c(\xi)}{|\xi|^2}\nonumber\\
&-c(Z)A(Y)X_n\xi_n\frac{c(\xi)}{|\xi|^2}-c(Z)\sum_{l=1}^{n-1}A(X)Y_l\xi_l\frac{c(\xi)}{|\xi|^2}-c(Z)A(X)Y_n\xi_n\frac{c(\xi)}{|\xi|^2}.
\end{align}
If we omit some items that have no contribution for computing ${\bf \Phi_5}$, by the Cauchy integral formula, we obtain
\begin{align}\label{64}
\pi^+_{\xi_n}A^1(x_0)=&-\sum_{j=1}^nX_j\frac{\partial {Y_n}}{\partial {x_j}}\bigg[\frac{c(Z)c(\xi')}{2(\xi_n-i)}+\frac{ic(Z)c(dx_n)}{2(\xi_n-i)}\bigg]\nonumber\\
&-X_n\bigg[\frac{c(Z)A(Y)c(\xi')}{2(\xi_n-i)}+\frac{ic(Z)A(Y)c(dx_n)}{2(\xi_n-i)}\bigg]\nonumber\\
&-Y_n\bigg[\frac{c(Z)A(X)c(\xi')}{2(\xi_n-i)}+\frac{ic(Z)A(X)c(dx_n)}{2(\xi_n-i)}\bigg].
\end{align}
By the relation of the Clifford action and ${\rm tr}{(AB)}={\rm tr }{(BA)}$,  we have the equalities:
\begin{align}\label{a49111}
&{\rm tr}[c(Z)A(X)c(\xi')]=\frac{1}{2}\sum_{j=1}^{n-1} \langle\nabla_X^L\frac{\partial}{\partial {x_n}},e_j\rangle \xi_j{\rm tr}[\texttt{id}];~~{\rm tr}[c(Z)A(X)c(dx_n)]=0;\nonumber\\
&{\rm tr}[c(Z)A(Y)c(\xi')]=\frac{1}{2}\sum_{j=1}^{n-1} \langle\nabla_Y^L\frac{\partial}{\partial {x_n}},e_j\rangle \xi_j{\rm tr}[\texttt{id}];~~{\rm
tr}[c(Z)A(Y)c(dx_n)]=0.
\end{align}
We note that $i<n,~\int_{|\xi'|=1}\{\xi_{i_{1}}\xi_{i_{2}}\cdots\xi_{i_{2q+1}}\}\sigma(\xi')=0$. Then, we have
\begin{equation}
{\rm tr} [\pi^+_{\xi_n}A^1\times
\partial_{\xi_n}\sigma_{-2m}(D^{-2m})](x_0)|_{|\xi'|=1}=-\sum_{j=1}^nX_j\frac{\partial {Y_n}}{\partial {x_j}}\frac{mi\xi_n}{(\xi_n-i)^{m+2}(\xi_n+i)^{m+1}}Z_n{\rm tr}[\texttt{id}].
\end{equation}
Then
\begin{align}\label{65}
&-i\int_{|\xi'|=1}\int^{+\infty}_{-\infty}{\rm tr} [\pi^+_{\xi_n}A^1\times
\partial_{\xi_n}\sigma_{-2m}(D^{-2m})](x_0)d\xi_n\sigma(\xi')dx'\nonumber\\
=&i\int_{|\xi'|=1}\int^{+\infty}_{-\infty}\sum_{j=1}^nX_j\frac{\partial {Y_n}}{\partial {x_j}}\frac{mi\xi_n}{(\xi_n-i)^{m+2}(\xi_n+i)^{m+1}}Z_n{\rm tr}[\texttt{id}]d\xi_n\sigma(\xi')dx'\nonumber\\
=&-Vol(S^{n-2})X(Y_n)Z_n2^m\int_{\Gamma^{+}}\frac{m\xi_n}{(\xi_n-i)^{m+2}(\xi_n+i)^{m+1}}d\xi_ndx'\nonumber\\
=&-Vol(S^{n-2})X(Y_n)Z_n2^m\frac{\pi}{(m+1)!}\left[\frac{2mi\xi_n}{(\xi_n+i)^{m+1}}\right]^{(m+1)}\bigg|_{\xi_n=i}dx'\nonumber\\
:=&-X(Y_n)Z_n\frac{\pi}{(m+1)!}D_0Vol(S^{n-2})2^mdx',
\end{align}
where
\begin{equation}
D_0=\left[\frac{2mi\xi_n}{(\xi_n+i)^{m+1}}\right]^{(m+1)}\bigg|_{\xi_n=i}=\frac{m(m+1)\cdot\cdot\cdot(2m)}{2^{2m+1}}.
\end{equation}
We note that
\begin{align}\label{45}
\sigma_{0}(D)(x_0)&=-\frac{1}{4}\sum_{s,t,i}\omega_{s,t}(e_i)
(x_{0})c(e_i)c(e_s)c(e_t):=H(x_0).
\end{align}
Secondly, by $H(x_0)=c_0c(dx_n)=-\frac{m}{2}h'(0)c(dx_n)$, further calculation leads to new results
\begin{align}
A^2(x_0)=&-\sum_{j,l=1}^{n}X_jY_l\xi_j\xi_l\bigg[\frac{[m\xi_n^4+2\xi_n^2-(m+2)]h'(0)c(Z)c(dx_n)}{2(1+\xi_n^2)^3}+\frac{[m\xi_n^3+(m+2)\xi_n]h'(0)c(Z)c(\xi')}{(1+\xi_n^2)^3}\nonumber\\
&+\frac{c(Z)c(\xi')c(dx_n)\partial_{x_n}[c(\xi')]}{(1+\xi_n^2)^2}-\frac{\xi_nc(Z)\partial_{x_n}[c(\xi')]}{(1+\xi_n^2)^2}\bigg].
\end{align}
If we omit some items that have no contribution for computing ${\bf \Phi_5}$, by the Cauchy integral formula, we obtain
\begin{align}\label{64}
\pi^+_{\xi_n}A^2(x_0)=&-\sum_{j,l=1}^{n-1}X_jY_l\xi_j\xi_l\bigg[\frac{[i\xi_n^2+(2m+3)\xi_n-(2m+4)i]h'(0)c(Z)c(dx_n)}{8(\xi_n-i)^3}+\frac{ic(Z)\partial_{x_n}[c(\xi')](x_0)}{4(\xi_n-i)^2}\nonumber\\
&-\frac{[(2m+1)i\xi_n+(2m+3)]h'(0)c(Z)c(\xi')}{8(\xi_n-i)^3}-\frac{(2+i\xi_n)c(Z)c(\xi')c(dx_n)\partial_{x_n}[c(\xi')](x_0)}{4(\xi_n-i)^2}\bigg]\nonumber\\
&-X_nY_n\bigg[\frac{[(4m-1)i\xi_n^2+(6m+1)\xi_n-2mi]h'(0)c(Z)c(dx_n)}{8(\xi_n-i)^3}-\frac{(2\xi_n-i)c(Z)\partial_{x_n}[c(\xi')](x_0)}{4(\xi_n-i)^2}\nonumber\\
&+\frac{[4m\xi_n^2-(6m+3)i\xi_n-(2m+1)]h'(0)c(Z)c(\xi')}{8(\xi_n-i)^3}-\frac{i\xi_nc(Z)c(\xi')c(dx_n)\partial_{x_n}[c(\xi')](x_0)}{4(\xi_n-i)^2}\bigg].\nonumber\\
\end{align}
Because ${\rm tr}[c(Z)c(\xi')c(dx_n)\partial_{x_n}[c(\xi')]]=-\frac{h'(0)}{2}Z_n{\rm tr}[\texttt{id}]$, and $g(Z,\xi')$ has no contribution for computing ${\bf \Phi_5}$. Then, we obtain
\begin{align}\label{ggg32}
&{\rm tr} [\pi^+_{\xi_n}A^2\times
\partial_{\xi_n}\sigma_{-2m}(D^{-2m})](x_0)|_{|\xi'|=1}\nonumber\\
=&-\sum_{j,l=1}^{n-1}X_jY_l\xi_j\xi_l\frac{m^2\xi_n^2-(m^2+m)i\xi_n}{2(\xi_n-i)^{m+4}(\xi_n+i)^{m+1}}Z_nh'(0){\rm tr}[\texttt{id}]\nonumber\\
&-X_nY_n\frac{(2m^2-m)i\xi_n^3+3m^2\xi_n^2-m^2i\xi_n}{2(\xi_n-i)^{m+4}(\xi_n+i)^{m+1}}Z_nh'(0){\rm tr}[\texttt{id}].
\end{align}
Then
\begin{align}\label{65}
&-i\int_{|\xi'|=1}\int^{+\infty}_{-\infty}{\rm tr} [\pi^+_{\xi_n}A^2\times
\partial_{\xi_n}\sigma_{-2m}(D^{-2m})](x_0)d\xi_n\sigma(\xi')dx'\nonumber\\
=&i\int_{|\xi'|=1}\int^{+\infty}_{-\infty}\sum_{j,l=1}^{n-1}X_jY_l\xi_j\xi_l\frac{m^2\xi_n^2-(m^2+m)i\xi_n}{2(\xi_n-i)^{m+4}(\xi_n+i)^{m+1}}Z_nh'(0){\rm tr}[\texttt{id}]d\xi_n\sigma(\xi')dx'\nonumber\\
&+i\int_{|\xi'|=1}\int^{+\infty}_{-\infty}X_nY_n\frac{(2m^2-m)i\xi_n^3+3m^2\xi_n^2-m^2i\xi_n}{2(\xi_n-i)^{m+4}(\xi_n+i)^{m+1}}Z_nh'(0){\rm tr}[\texttt{id}]d\xi_n\sigma(\xi')dx'\nonumber\\
=&\frac{i}{4}Vol(S^{n-2})g(X^T,Y^T)Z_nh'(0)2^{m}\int_{\Gamma^{+}}\frac{m\xi_n^2-(m+1)i\xi_n}{(\xi_n-i)^{m+4}(\xi_n+i)^{m+1}}d\xi_ndx'\nonumber\\
&+\frac{i}{2}Vol(S^{n-2})X_nY_nZ_nh'(0)2^{m}\int_{\Gamma^{+}}\frac{(2m^2-m)i\xi_n^3+3m^2\xi_n^2-m^2i\xi_n}{(\xi_n-i)^{m+4}(\xi_n+i)^{m+1}}d\xi_ndx'\nonumber\\
=&\frac{i}{4}Vol(S^{n-2})g(X^T,Y^T)Z_nh'(0)2^{m}\frac{2\pi i}{(m+3)!}\left[\frac{m\xi_n^2-(m+1)i\xi_n}{(\xi_n+i)^{m+1}}\right]^{(m+3)}\bigg|_{\xi_n=i}dx'\nonumber\\
&+\frac{i}{2}Vol(S^{n-2})X_nY_nZ_nh'(0)2^{m}\frac{2\pi i}{(m+3)!}\left[\frac{(2m^2-m)i\xi_n^3+3m^2\xi_n^2-m^2i\xi_n}{(\xi_n+i)^{m+1}}\right]^{(m+3)}\bigg|_{\xi_n=i}dx'\nonumber\\
:=&-\Big(g(X^T,Y^T)\frac{\pi}{2(m+3)!}D_1+X_nY_n\frac{\pi}{(m+3)!}D_2\Big)Vol(S^{n-2})Z_nh'(0)2^{m}dx',
\end{align}
where
\begin{align}
D_1=&\left[\frac{m\xi_n^2-(m+1)i\xi_n}{(\xi_n+i)^{m+1}}\right]^{(m+3)}\bigg|_{\xi_n=i}=\frac{(2m^2+9m+3)(m+1)(m+2)\cdot\cdot\cdot(2m+1)}{2^{2m+3}};\nonumber\\
D_2=&\left[\frac{(2m^2-m)i\xi_n^3+3m^2\xi_n^2-m^2i\xi_n}{(\xi_n+i)^{m+1}}\right]^{(m+3)}\bigg|_{\xi_n=i}=\frac{3(2m^2+9m-1)m(m+1)\cdot\cdot\cdot(2m)}{2^{2m+3}}.
\end{align}
Thirdly, for $A^3$, we get
\begin{align}\label{66pp}
\pi^+_{\xi_n} A^3(x_0)=&-X_nY_n\bigg[\frac{c(Z)\partial_{x_n}[c(\xi')](x_0)}{(\xi_n-i)}+h'(0)\bigg(\frac{ic(Z)c(\xi')}{2(\xi_n-i)^2}+\frac{i\xi_nc(Z)c(dx_n)}{2(\xi_n-i)^2}\bigg)\bigg].
\end{align}
Moreover
\begin{align}
{\rm tr} [\pi^+_{\xi_n}A^3\times
\partial_{\xi_n}\sigma_{-2m}(D^{-2m})](x_0)|_{|\xi'|=1}=-X_nY_n\frac{mi\xi_n^2}{(\xi_n-i)^{m+3}(\xi_n+i)^{m+1}}Z_nh'(0){\rm tr}[\texttt{id}].
\end{align}
Then
\begin{align}
&-i\int_{|\xi'|=1}\int^{+\infty}_{-\infty}{\rm tr} [\pi^+_{\xi_n}A^3\times
\partial_{\xi_n}\sigma_{-2m}(D^{-2m})](x_0)d\xi_n\sigma(\xi')dx'\nonumber\\
=&i\int_{|\xi'|=1}\int^{+\infty}_{-\infty}X_nY_n\frac{mi\xi_n^2}{(\xi_n-i)^{m+3}(\xi_n+i)^{m+1}}Z_nh'(0){\rm tr}[\texttt{id}]d\xi_n\sigma(\xi')dx'\nonumber\\
=&-Vol(S^{n-2})X_nY_nZ_nh'(0)2^m\int_{\Gamma^+}\frac{m\xi_n^2}{(\xi_n-i)^{m+3}(\xi_n+i)^{m+1}}d\xi_ndx'\nonumber\\
=&-Vol(S^{n-2})X_nY_nZ_nh'(0)2^m\frac{\pi}{(m+2)!}\left[\frac{2mi\xi_n^2}{(\xi_n+i)^{m+1}}\right]^{(m+2)}\bigg|_{\xi_n=i}dx'\nonumber\\
:=&-X_nY_n\frac{\pi}{(m+2)!}D_3Vol(S^{n-2})Z_nh'(0)2^mdx',
\end{align}
where
\begin{align}
D_3=\left[\frac{2mi\xi_n^2}{(\xi_n+i)^{m+1}}\right]^{(m+2)}\bigg|_{\xi_n=i}=-\frac{(m-1)m\cdot\cdot\cdot(2m)}{2^{2m+1}}.
\end{align}
Therefore, we obtain
\begin{align}\label{6666}
\Phi_5=&-i\int_{|\xi'|=1}\int^{+\infty}_{-\infty}{\rm tr} [\pi^+_{\xi_n}(A^1+A^2+A^3)\times
\partial_{\xi_n}\sigma_{-2m}(D^{-2m})](x_0)d\xi_n\sigma(\xi')dx'\nonumber\\
=&\bigg\{X(Y_n)Z_n\frac{\pi}{(m+1)!}D_0+Z_nh'(0)\Big(g(X^T,Y^T)\frac{\pi}{2(m+3)!}D_1+X_nY_n\frac{\pi}{(m+3)!}D_2\nonumber\\
&+X_nY_n\frac{\pi}{(m+2)!}D_3\Big)\bigg\}Vol(S^{n-2})2^mdx'.
\end{align}
Now $\Phi$ is the sum of the the {\bf  Case (a-I)}-{\bf  Case (a-V)}. Therefore, we get
\begin{align}\label{795}
\Phi&=\sum_{i=1}^5\Phi_i\nonumber\\
=&\bigg\{\partial_{x_n}\Big(g(X^T,Y^T)Z_n\Big)\frac{\pi}{2(m+2)!}A_0-\partial_{x_n}\Big((X_nY_n)Z_n\Big)\frac{\pi}{(m+2)!}A_1+g(X^T,Y^T)Z_nh'(0)\Big(\frac{\pi i}{4(m+3)!}A_2\nonumber\\
&+\frac{\pi}{2(m+3)!}B_0+\frac{\pi}{4(m+2)!}C_0+\frac{\pi}{2(m+3)!}C_1-\frac{\pi}{2(m+3)!}D_1\Big)-X_nY_nZ_nh'(0)\Big(-\frac{\pi i}{2(m+3)!}A_3\nonumber\\
&+\frac{\pi}{(m+3)!}B_1+\frac{\pi}{2(m+2)!}C_2+\frac{\pi}{(m+3)!}C_3+\frac{\pi}{(m+3)!}D_2+\frac{\pi}{(m+2)!}D_3\Big)\nonumber\\
&-X(Y_n)Z_n\frac{\pi}{(m+1)!}D_0\bigg\}Vol(S^{n-2})2^mdx'.
\end{align}
Then, by (\ref{b141}) and (\ref{795}), we obtain the following theorem.
\begin{thm}\label{thmb1}
Let $M$ be an $n=(2m+1)$-dimensional oriented
compact spin manifold with the boundary $\partial M$, then we get the following equality:
\begin{align}
\label{b2773}
&\widetilde{{\rm Wres}}[\pi^+(c(Z)\nabla_X^{S(TM)}\nabla_Y^{S(TM)}D^{-1})\circ\pi^+(D^{-2m})]\nonumber\\
=&\int_{\partial M}\bigg\{\partial_{x_n}\Big(g(X^T,Y^T)Z_n\Big)\frac{\pi}{2(m+2)!}A_0-\partial_{x_n}\Big((X_nY_n)Z_n\Big)\frac{\pi}{(m+2)!}A_1+g(X^T,Y^T)Z_nh'(0)\Big(\frac{\pi i}{4(m+3)!}A_2\nonumber\\
&+\frac{\pi}{2(m+3)!}B_0+\frac{\pi}{4(m+2)!}C_0+\frac{\pi}{2(m+3)!}C_1-\frac{\pi}{2(m+3)!}D_1\Big)-X_nY_nZ_nh'(0)\Big(-\frac{\pi i}{2(m+3)!}A_3\nonumber\\
&+\frac{\pi}{(m+3)!}B_1+\frac{\pi}{2(m+2)!}C_2+\frac{\pi}{(m+3)!}C_3+\frac{\pi}{(m+3)!}D_2+\frac{\pi}{(m+2)!}D_3\Big)\nonumber\\
&-X(Y_n)Z_n\frac{\pi}{(m+1)!}D_0\bigg\}Vol(S^{n-2})2^md{\rm Vol_{M}}.
\end{align}
\end{thm}

\section{The noncommutative residue $\widetilde{{\rm Wres}}[\pi^+(c(Z)\nabla_X^{S(TM)}\nabla_Y^{S(TM)}D^{-2})\circ\pi^+(D^{-(2m-1)})]$ on odd dimensional manifolds with boundary}
\label{section:4}
Next, we compute the noncommutative residue $\widetilde{{\rm Wres}}[\pi^+(c(Z)\nabla_X^{S(TM)}\nabla_Y^{S(TM)}D^{-2})\circ\pi^+(D^{-(2m-1)})]$ on $(2m+1)$-dimensional oriented compact spin manifolds with boundary and get a general Dabrowski-Sitarz-Zalecki type theorem in this case.

Similar to \cite{Wa3}, by $\int_{S(\xi)}\xi_{j_{1}}\cdots\xi_{j_{l}}\texttt{d}(\xi)=0$ for odd $l$, we can compute the noncommutative residue
\begin{align}
\label{c121}
\widetilde{{\rm Wres}}[\pi^+(c(Z)\nabla_X^{S(TM)}\nabla_Y^{S(TM)}D^{-2})\circ\pi^+(D^{-(2m-1)})]=\int_{\partial M}\widetilde{\Phi},
\end{align}
where
\begin{align}
\label{c2}
 \widetilde{\Phi}=&\int_{|\xi'|=1}\int^{+\infty}_{-\infty}\sum^{\infty}_{j, k=0}\sum\frac{(-i)^{|\alpha|+j+k+1}}{\alpha!(j+k+1)!}
\times {\rm trace}_{S(TM)}[\partial^j_{x_n}\partial^\alpha_{\xi'}\partial^k_{\xi_n}\sigma^+_{r}(c(Z)\nabla_X^{S(TM)}\nabla_Y^{S(TM)}D^{-2})(x',0,\xi',\xi_n)\nonumber\\
&\times\partial^\alpha_{x'}\partial^{j+1}_{\xi_n}\partial^k_{x_n}\sigma_{l}(D^{-(2m-1)})(x',0,\xi',\xi_n)]d\xi_n\sigma(\xi')dx',
\end{align}
and the sum is taken over $r+l-k-j-|\alpha|-1=-(2m+1),~~r\leq 0,~~l\leq -(2m-1)$.

Next, we need to compute $\int_{\partial M} \widetilde{\Phi}$. When sum is taken over $
r+l-k-j-|\alpha|=-2m,~~r\leq 0,~~l\leq -(2m-1),$ then we have the following five cases:\\

\noindent  {\bf Case (b-I)}~$r=0,~l=-(2m-1),~k=j=0,~|\alpha|=1$.\\

\noindent By (\ref{c2}), we get
\begin{equation}
\label{c4}
\widetilde{\Phi}_1=-\int_{|\xi'|=1}\int^{+\infty}_{-\infty}\sum_{|\alpha|=1}
{\rm tr}[\partial^\alpha_{\xi'}\pi^+_{\xi_n}\sigma_{0}(c(Z)\nabla_X^{S(TM)}\nabla_Y^{S(TM)}D^{-2})\times
 \partial^\alpha_{x'}\partial_{\xi_n}\sigma_{-(2m-1)}(D^{-(2m-1)})](x_0)d\xi_n\sigma(\xi')dx'.
\end{equation}
Similarly,
\begin{align}
\label{c5}
\partial_{x_i}\sigma_{-(2m-1)}(D^{-(2m-1)})(x_0)&=\partial_{x_i}{(ic(\xi)|\xi|^{-2m})}(x_0)\nonumber\\
&=i\partial_{x_i}c(\xi)(x_0)|\xi|^{-2m}+ic(\xi)\partial_{x_i}(|\xi|^{-2m})(x_0)\nonumber\\
 &=0,
\end{align}
\noindent so $\widetilde{\Phi}_1=0$.\\

 \noindent  {\bf Case (b-II)}~$r=0,~l=-(2m-1),~k=|\alpha|=0,~j=1$.\\

\noindent By (\ref{c2}), we get
\begin{align}
\label{c6}
\widetilde{\Phi}_2&=-\frac{1}{2}\int_{|\xi'|=1}\int^{+\infty}_{-\infty} {\rm
tr} [\partial_{x_n}\pi^+_{\xi_n}\sigma_{0}(c(Z)\nabla_X^{S(TM)}\nabla_Y^{S(TM)}D^{-2})\times
\partial_{\xi_n}^2\sigma_{-(2m-1)}(D^{-(2m-1)})](x_0)d\xi_n\sigma(\xi')dx'\nonumber\\
&=-\frac{1}{2}\int_{|\xi'|=1}\int^{+\infty}_{-\infty} {\rm
tr} [\partial_{\xi_n}^2\partial_{x_n}\pi^+_{\xi_n}\sigma_{0}(c(Z)\nabla_X^{S(TM)}\nabla_Y^{S(TM)}D^{-2})\times
\sigma_{-(2m-1)}(D^{-(2m-1)})](x_0)d\xi_n\sigma(\xi')dx'.
\end{align}
\noindent By (3.11) in \cite{Wa7}, we have\\
\begin{eqnarray}\label{c7}
\sigma_{-(2m-1)}(D^{-(2m-1)})(x_0)=\frac{i[c(\xi')+\xi_nc(\mathrm{d}x_n)]}{(1+\xi_n^2)^m}.
\end{eqnarray}
By Lemma \ref{lema}, we get
\begin{align}\label{c8}
&\partial_{x_n}\sigma_{0}(c(Z)\nabla_X^{S(TM)}\nabla_Y^{S(TM)}D^{-2})(x_0)\nonumber\\
=&\partial_{x_n}\bigg(-i\sum_{j,l=1}^nX_jY_l\xi_j\xi_l\frac{c(Z)}{|\xi|^{2}}\bigg)(x_0)\nonumber\\
=&-\sum_{j,l=1}^n\left(\frac{\partial{X_j}}{\partial {x_n}}Y_l+X_j\frac{\partial{Y_l}}{\partial {x_n}}\right)\xi_j\xi_l\frac{c(Z)}{|\xi|^{2}}-\sum_{j,l=1}^nX_jY_l\xi_j\xi_l \left[\frac{\partial_{x_n}(c(Z))}{|\xi|^2}-\frac{c(Z)|\xi'|^2h'(0)}{|\xi|^4}\right].
\end{align}
We note that $i<n,~\int_{|\xi'|=1}\{\xi_{i_{1}}\xi_{i_{2}}\cdots\xi_{i_{2q+1}}\}\sigma(\xi')=0$,
so we omit some items that have no contribution for computing ${\bf \widetilde{\Phi}_2}$. Then, we have
\begin{align}\label{c9}
&\partial_{x_n}\pi^+_{\xi_n}\sigma_{0}(c(Z)\nabla_X^{S(TM)}\nabla_Y^{S(TM)}D^{-2})(x_0)\nonumber\\
=&\pi^+_{\xi_n}\partial_{x_n}\sigma_{0}(c(Z)\nabla_X^{S(TM)}\nabla_Y^{S(TM)}D^{-2})(x_0)\nonumber\\
=&\bigg(\sum_{j,l=1}^{n-1}\frac{\partial{X_j}}{\partial {x_n}}Y_l\xi_j\xi_l-\frac{\partial{X_n}}{\partial {x_n}}Y_n+\sum_{j,l=1}^{n-1}X_j\frac{\partial{Y_l}}{\partial {x_n}}\xi_j\xi_l-X_n\frac{\partial{Y_n}}{\partial {x_n}}\bigg)\frac{ic(Z)}{2(\xi_n-i)}\nonumber\\
&+\sum_{j,l=1}^{n-1}X_jY_l\xi_j\xi_l\bigg(\frac{i\partial_{x_n}(c(Z))}{2(\xi_n-i)}-h'(0)|\xi'|^2\frac{(i\xi_n+2)c(Z)}{4(\xi_n-i)^2}\bigg)\nonumber\\
&-X_nY_n\bigg(\frac{i\partial_{x_n}(c(Z))}{2(\xi_n-i)}+h'(0)|\xi'|^2\frac{i\xi_nc(Z)}{4(\xi_n-i)^2}\bigg).
\end{align}
By further calculation, we have
\begin{align}\label{555}
&\partial_{\xi_n}^2\partial_{x_n}\pi^+_{\xi_n}\sigma_{0}(c(Z)\nabla_X^{S(TM)}\nabla_Y^{S(TM)}D^{-2})(x_0)\nonumber\\
=&\bigg(\sum_{j,l=1}^{n-1}\frac{\partial{X_j}}{\partial {x_n}}Y_l\xi_j\xi_l-\frac{\partial{X_n}}{\partial {x_n}}Y_n+\sum_{j,l=1}^{n-1}X_j\frac{\partial{Y_l}}{\partial {x_n}}\xi_j\xi_l-X_n\frac{\partial{Y_n}}{\partial {x_n}}\bigg)\frac{ic(Z)}{(\xi_n-i)^3}\nonumber\\
&+\sum_{j,l=1}^{n-1}X_jY_l\xi_j\xi_l\bigg(\frac{i\partial_{x_n}(c(Z))}{(\xi_n-i)^3}-h'(0)|\xi'|^2\frac{(i\xi_n+4)c(Z)}{2(\xi_n-i)^4}\bigg)\nonumber\\
&-X_nY_n\bigg(\frac{i\partial_{x_n}(c(Z))}{(\xi_n-i)^3}+h'(0)|\xi'|^2\frac{(i\xi_n-2)c(Z)}{2(\xi_n-i)^4}\bigg).
\end{align}
Moreover
\begin{align}\label{35}
&{\rm tr} [\partial_{\xi_n}^2\partial_{x_n}\pi^+_{\xi_n}\sigma_{0}(c(Z)\nabla_X^{S(TM)}\nabla_Y^{S(TM)}D^{-2})\times
\sigma_{-(2m-1)}(D^{-(2m-1)})](x_0)|_{|\xi'|=1}\nonumber\\
=&\sum_{j,l=1}^{n-1}\bigg(\frac{\partial{X_j}}{\partial {x_n}}Y_l+X_j\frac{\partial{Y_l}}{\partial {x_n}}\bigg)\xi_j\xi_l\frac{\xi_n}{(\xi_n-i)^{m+3}(\xi_n+i)^{m}}Z_n{\rm tr}[\texttt{id}]\nonumber\\
&-\bigg(\frac{\partial{X_n}}{\partial {x_n}}Y_n+X_n\frac{\partial{Y_n}}{\partial {x_n}}\bigg)\frac{\xi_n}{(\xi_n-i)^{m+3}(\xi_n+i)^{m}}Z_n{\rm tr}[\texttt{id}]\nonumber\\
&+\sum_{j,l=1}^{n-1}X_jY_l\xi_j\xi_l\frac{\xi_n}{(\xi_n-i)^{m+3}(\xi_n+i)^{m}}\partial_{x_n}(Z_n){\rm tr}[\texttt{id}]\nonumber\\
&-\sum_{j,l=1}^{n-1}X_jY_l\xi_j\xi_l\frac{\xi_n^2-4i\xi_n}{2(\xi_n-i)^{m+4}(\xi_n+i)^{m}}Z_nh'(0){\rm tr}[\texttt{id}]\nonumber\\
&-X_nY_n\frac{\xi_n}{(\xi_n-i)^{m+3}(\xi_n+i)^{m}}\partial_{x_n}(Z_n){\rm tr}[\texttt{id}]\nonumber\\
&-X_nY_n\frac{\xi_n^2+2i\xi_n}{2(\xi_n-i)^{m+4}(\xi_n+i)^{m}}Z_nh'(0){\rm tr}[\texttt{id}].
\end{align}
Therefore, we get
\begin{align}\label{35kkk}
\widetilde{\Phi}_2
=&-\frac{1}{2}\int_{|\xi'|=1}\int^{+\infty}_{-\infty} \sum_{j,l=1}^{n-1}\bigg(\frac{\partial{X_j}}{\partial {x_n}}Y_l+X_j\frac{\partial{Y_l}}{\partial {x_n}}\bigg)\xi_j\xi_l\frac{\xi_n}{(\xi_n-i)^{m+3}(\xi_n+i)^{m}}Z_n{\rm tr}[\texttt{id}]d\xi_n\sigma(\xi')dx'\nonumber\\
&+\frac{1}{2}\int_{|\xi'|=1}\int^{+\infty}_{-\infty} \bigg(\frac{\partial{X_n}}{\partial {x_n}}Y_n+X_n\frac{\partial{Y_n}}{\partial {x_n}}\bigg)\frac{\xi_n}{(\xi_n-i)^{m+3}(\xi_n+i)^{m}}Z_n{\rm tr}[\texttt{id}]d\xi_n\sigma(\xi')dx'\nonumber\\
&-\frac{1}{2}\int_{|\xi'|=1}\int^{+\infty}_{-\infty} \sum_{j,l=1}^{n-1}X_jY_l\xi_j\xi_l\frac{\xi_n}{(\xi_n-i)^{m+3}(\xi_n+i)^{m}}\partial_{x_n}(Z_n){\rm tr}[\texttt{id}]d\xi_n\sigma(\xi')dx'\nonumber\\
&+\frac{1}{2}\int_{|\xi'|=1}\int^{+\infty}_{-\infty} \sum_{j,l=1}^{n-1}X_jY_l\xi_j\xi_l\frac{\xi_n^2-4i\xi_n}{2(\xi_n-i)^{m+4}(\xi_n+i)^{m}}Z_nh'(0){\rm tr}[\texttt{id}]d\xi_n\sigma(\xi')dx'\nonumber\\
&+\frac{1}{2}\int_{|\xi'|=1}\int^{+\infty}_{-\infty} X_nY_n\frac{\xi_n}{(\xi_n-i)^{m+3}(\xi_n+i)^{m}}\partial_{x_n}(Z_n){\rm tr}[\texttt{id}]d\xi_n\sigma(\xi')dx'\nonumber\\
&+\frac{1}{2}\int_{|\xi'|=1}\int^{+\infty}_{-\infty} X_nY_n\frac{\xi_n^2+2i\xi_n}{2(\xi_n-i)^{m+4}(\xi_n+i)^{m}}Z_nh'(0){\rm tr}[\texttt{id}]d\xi_n\sigma(\xi')dx'\nonumber\\
=&-\frac{1}{4m}Vol(S^{n-2})\partial_{x_n}\Big(g(X^T,Y^T)Z_n\Big)2^m\int_{\Gamma^{+}}\frac{\xi_n}{(\xi_n-i)^{m+3}(\xi_n+i)^{m}}d\xi_ndx'\nonumber\\
&+\frac{1}{2}Vol(S^{n-2})\partial_{x_n}\Big((X_nY_n)Z_n\Big)2^m\int_{\Gamma^{+}}\frac{\xi_n}{(\xi_n-i)^{m+3}(\xi_n+i)^{m}}d\xi_ndx'\nonumber\\
&+\frac{1}{8m}Vol(S^{n-2})g(X^T,Y^T)Z_nh'(0)2^m\int_{\Gamma^{+}}\frac{\xi_n^2-4i\xi_n}{(\xi_n-i)^{m+4}(\xi_n+i)^{m}}d\xi_ndx'\nonumber\\
&+\frac{1}{4}Vol(S^{n-2})X_nY_nZ_nh'(0)2^m\int_{\Gamma^{+}}\frac{\xi_n^2+2i\xi_n}{(\xi_n-i)^{m+4}(\xi_n+i)^{m}}d\xi_ndx'\nonumber\\
=&-\frac{1}{4}Vol(S^{n-2})\partial_{x_n}\Big(g(X^T,Y^T)Z_n\Big)2^m\frac{2\pi i}{(m+2)!}\left[\frac{\xi_n}{m(\xi_n+i)^{m}}\right]^{(m+2)}\bigg|_{\xi_n=i}dx'\nonumber\\
&+\frac{1}{2}Vol(S^{n-2})\partial_{x_n}\Big((X_nY_n)Z_n\Big)2^m\frac{2\pi i}{(m+2)!}\left[\frac{\xi_n}{(\xi_n+i)^{m}}\right]^{(m+2)}\bigg|_{\xi_n=i}dx'\nonumber\\
&+\frac{1}{8}Vol(S^{n-2})g(X^T,Y^T)Z_nh'(0)2^m\frac{2\pi i}{(m+3)!}\left[\frac{\xi_n^2-4i\xi_n}{m(\xi_n+i)^{m}}\right]^{(m+3)}\bigg|_{\xi_n=i}dx'\nonumber\\
&+\frac{1}{4}Vol(S^{n-2})X_nY_nZ_nh'(0)2^m\frac{2\pi i}{(m+3)!}\left[\frac{\xi_n^2+2i\xi_n}{(\xi_n+i)^{m}}\right]^{(m+3)}\bigg|_{\xi_n=i}dx'\nonumber\\
:=&\bigg\{\partial_{x_n}\Big(g(X^T,Y^T)Z_n\Big)\frac{-\pi i}{2(m+2)!}E_0+\partial_{x_n}\Big((X_nY_n)Z_n\Big)\frac{\pi i}{(m+2)!}E_1\nonumber\\
&+Z_nh'(0)\Big(g(X^T,Y^T)\frac{\pi i}{4(m+3)!}E_2+X_nY_n\frac{\pi i}{2(m+3)!}E_3\Big)\bigg\}Vol(S^{n-2})2^mdx',
\end{align}
where
\begin{align}
E_0&=\left[\frac{\xi_n}{m(\xi_n+i)^{m}}\right]^{(m+2)}\bigg|_{\xi_n=i}=-\frac{3(m+1)(m+2)\cdot\cdot\cdot(2m)}{2^{2m+2}i};\nonumber\\
E_1&=\left[\frac{\xi_n}{(\xi_n+i)^{m}}\right]^{(m+2)}\bigg|_{\xi_n=i}=-\frac{3m(m+1)\cdot\cdot\cdot(2m)}{2^{2m+2}i};\nonumber\\
E_2&=\left[\frac{\xi_n^2-4i\xi_n}{m(\xi_n+i)^{m}}\right]^{(m+3)}\bigg|_{\xi_n=i}=-\frac{15(m+1)(m+1)(m+2)\cdot\cdot\cdot(2m)}{2^{2m+2}i};\nonumber\\
E_3&=\left[\frac{\xi_n^2+2i\xi_n}{(\xi_n+i)^{m}}\right]^{(m+3)}\bigg|_{\xi_n=i}=\frac{3(3m-1)m(m+1)\cdot\cdot\cdot(2m)}{2^{2m+2}i}.
\end{align}\\

\noindent  {\bf Case (b-III)}~$r=0,~l=-(2m-1),~j=|\alpha|=0,~k=1$.\\

\noindent By (\ref{c2}), we get
\begin{align}\label{c12}
\widetilde{\Phi}_3&=-\frac{1}{2}\int_{|\xi'|=1}\int^{+\infty}_{-\infty}
{\rm tr} [\partial_{\xi_n}\pi^+_{\xi_n}\sigma_{0}(c(Z)\nabla_X^{S(TM)}\nabla_Y^{S(TM)}D^{-2})\times
\partial_{\xi_n}\partial_{x_n}\sigma_{-(2m-1)}(D^{-(2m-1)})](x_0)d\xi_n\sigma(\xi')dx'\nonumber\\
&=\frac{1}{2}\int_{|\xi'|=1}\int^{+\infty}_{-\infty}
{\rm tr} [\partial_{\xi_n}^2\pi^+_{\xi_n}\sigma_{0}(c(Z)\nabla_X^{S(TM)}\nabla_Y^{S(TM)}D^{-2})\times
\partial_{x_n}\sigma_{-(2m-1)}(D^{-(2m-1)})](x_0)d\xi_n\sigma(\xi')dx'.
\end{align}
\noindent By (3.17) in \cite{Wa7}, we have
\begin{eqnarray}\label{c13}
\partial_{x_n} \sigma_{-(2m-1)}(D^{-(2m-1)})(x_0)
=\frac{i\partial_{x_n}[c(\xi')](x_0)}{(1+\xi_n^2)^{m}}-\frac{mih'(0)c(\xi)}{(1+\xi_n^2)^{m+1}}.
\end{eqnarray}
By Lemma \ref{lema}, we have
\begin{align}\label{aa38}
\pi^+_{\xi_n}\sigma_{0}(c(Z)\nabla_X^{S(TM)}\nabla_Y^{S(TM)}D^{-2})(x_0)=&\bigg(\sum_{j,l=1}^{n-1}X_jY_l\xi_j\xi_l-X_nY_n\bigg)\frac{ic(Z)}{2(\xi_n-i)}\nonumber\\
&-\sum_{j=1}^{n-1}(X_jY_n+X_nY_j)\xi_j\frac{c(Z)}{2(\xi_n-i)}.
\end{align}
We omit some items that have no contribution for computing ${\bf \widetilde{\Phi}_3}$. Then, we get
\begin{align}\label{mmmmm}
\partial_{\xi_n}^2\pi^+_{\xi_n}\sigma_{0}(c(Z)\nabla_X^{S(TM)}\nabla_Y^{S(TM)}D^{-2})(x_0)=\bigg(\sum_{j,l=1}^{n-1}X_jY_l\xi_j\xi_l-X_nY_n\bigg)\frac{ic(Z)}{(\xi_n-i)^3}.
\end{align}
Moreover
\begin{align}\label{35}
&{\rm tr} [\partial_{\xi_n}^2\pi^+_{\xi_n}\sigma_{0}(c(Z)\nabla_X^{S(TM)}\nabla_Y^{S(TM)}D^{-2})\times
\partial_{x_n}\sigma_{-(2m-1)}(D^{-(2m-1)})](x_0)|_{|\xi'|=1}\nonumber\\
=&-\sum_{j,l=1}^{n-1}X_jY_l\xi_j\xi_l\frac{m\xi_n}{(\xi_n-i)^{m+4}(\xi_n+i)^{m+1}}Z_nh'(0){\rm tr}[\texttt{id}]\nonumber\\
&+X_nY_n\frac{m\xi_n}{(\xi_n-i)^{m+4}(\xi_n+i)^{m+1}}Z_nh'(0){\rm tr}[\texttt{id}].
\end{align}
Therefore, we get
\begin{align}\label{41}
\widetilde{\Phi}_3=&-\frac{1}{2}\int_{|\xi'|=1}\int^{+\infty}_{-\infty} \sum_{j,l=1}^{n-1}X_jY_l\xi_j\xi_l\frac{m\xi_n}{(\xi_n-i)^{m+4}(\xi_n+i)^{m+1}}Z_nh'(0){\rm tr}[\texttt{id}]d\xi_n\sigma(\xi')dx'\nonumber\\
&+\frac{1}{2}\int_{|\xi'|=1}\int^{+\infty}_{-\infty} X_nY_n\frac{m\xi_n}{(\xi_n-i)^{m+4}(\xi_n+i)^{m+1}}Z_nh'(0){\rm tr}[\texttt{id}]d\xi_n\sigma(\xi')dx'\nonumber\\
=&-\frac{1}{4}Vol(S^{n-2})g(X^T,Y^T)Z_nh'(0)2^m\int_{\Gamma^{+}} \frac{\xi_n}{(\xi_n-i)^{m+4}(\xi_n+i)^{m+1}}d\xi_ndx'\nonumber\\
&+\frac{1}{2}Vol(S^{n-2})X_nY_nZ_nh'(0)2^m\int_{\Gamma^{+}} \frac{m\xi_n}{(\xi_n-i)^{m+4}(\xi_n+i)^{m+1}}d\xi_ndx'\nonumber\\
=&-\frac{1}{4}Vol(S^{n-2})g(X^T,Y^T)Z_nh'(0)2^m\frac{2\pi i}{(m+3)!}\left[\frac{\xi_n}{(\xi_n+i)^{m+1}}\right]^{(m+3)}\bigg|_{\xi_n=i}dx'\nonumber\\
&+\frac{1}{2}Vol(S^{n-2})X_nY_nZ_nh'(0)2^m\frac{2\pi i}{(m+3)!}\left[\frac{m\xi_n}{(\xi_n+i)^{m+1}}\right]^{(m+3)}\bigg|_{\xi_n=i}dx'\nonumber\\
:=&\Big(-g(X^T,Y^T)\frac{\pi i}{2(m+3)!}F_0+X_nY_n\frac{\pi i}{(m+3)!}F_1\Big)Vol(S^{n-2})Z_nh'(0)2^mdx',
\end{align}
where
\begin{align}
F_0&=\left[\frac{\xi_n}{(\xi_n+i)^{m+1}}\right]^{(m+3)}\bigg|_{\xi_n=i}=-\frac{3(m+1)(m+2)\cdot\cdot\cdot(2m+2)}{2^{2m+4}i};\nonumber\\
F_1&=\left[\frac{m\xi_n}{(\xi_n+i)^{m+1}}\right]^{(m+3)}\bigg|_{\xi_n=i}=-\frac{3m(m+1)\cdot\cdot\cdot(2m+2)}{2^{2m+4}i}.
\end{align}\\

\noindent  {\bf Case (b-IV)}~$r=0,~l=-2m,~k=j=|\alpha|=0$.\\

\noindent By (\ref{c2}), we get
\begin{align}\label{c19}
\widetilde{\Phi}_4=&-i\int_{|\xi'|=1}\int^{+\infty}_{-\infty}{\rm tr} [\pi^+_{\xi_n}\sigma_{0}(c(Z)\nabla_X^{S(TM)}\nabla_Y^{S(TM)}D^{-2})\times
\partial_{\xi_n}\sigma_{-2m}(D^{-(2m-1)})](x_0)d\xi_n\sigma(\xi')dx'\nonumber\\
=&i\int_{|\xi'|=1}\int^{+\infty}_{-\infty}{\rm tr} [\partial_{\xi_n}\pi^+_{\xi_n}\sigma_{0}(c(Z)\nabla_X^{S(TM)}\nabla_Y^{S(TM)}D^{-2})\times
\sigma_{-2m}(D^{-(2m-1)})](x_0)d\xi_n\sigma(\xi')dx'.
\end{align}
 By (3.37) in \cite{Wa7}, we have
\begin{align}\label{c20}
\sigma_{-2m}(D^{-(2m-1)})(x_0)=&-\frac{(2m+1)h'(0)c(dx_n)}{4(1+\xi_n^2)^{m}}-\frac{2m\xi_n\partial_{x_n}[c(\xi')](x_0)}{(1+\xi_n^2)^{m+1}}\nonumber\\
&+\frac{m}{(1+\xi_n^2)^{m+1}}\bigg[\frac{1}{2}h'(0)c(\xi')c(dx_n)+\frac{2m+1}{2}h'(0)\xi_n\bigg]c(\xi)+\frac{(m^2+m)h'(0)\xi_n}{(1+\xi_n^2)^{m+2}}c(\xi).
\end{align}
By (\ref{aa38}), we have
\begin{align}\label{c21}
\partial_{\xi_n}\pi^+_{\xi_n}\sigma_{0}(c(Z)\nabla_X^{S(TM)}\nabla_Y^{S(TM)}D^{-2})(x_0)=&-\bigg(\sum_{j,l=1}^{n-1}X_jY_l\xi_j\xi_l-X_nY_n\bigg)\frac{ic(Z)}{2(\xi_n-i)^2}\nonumber\\
&+\sum_{j=1}^{n-1}(X_jY_n+X_nY_j)\xi_j\frac{c(Z)}{2(\xi_n-i)^2}.
\end{align}
We note that $i<n,~\int_{|\xi'|=1}\{\xi_{i_{1}}\xi_{i_{2}}\cdots\xi_{i_{2q+1}}\}\sigma(\xi')=0$, so we omit some items that have no contribution for computing ${\bf \widetilde{\Phi}_4}$. Then, we have
\begin{align}\label{35}
&{\rm tr} [\partial_{\xi_n}\pi^+_{\xi_n}\sigma_{0}(c(Z)\nabla_X^{S(TM)}\nabla_Y^{S(TM)}D^{-2})\times
\partial_{x_n}\sigma_{-2m}(D^{-(2m-1)})](x_0)|_{|\xi'|=1}\nonumber\\
=&-\sum_{j,l=1}^{n-1}X_jY_l\xi_j\xi_l\frac{(-4m^2+1)i\xi_n^4-(8m^2+4m-2)i\xi_n^2+i}{8(\xi_n-i)^{m+4}(\xi_n+i)^{m+2}}Z_nh'(0){\rm tr}[\texttt{id}]\nonumber\\
&+X_nY_n\frac{(-4m^2+1)i\xi_n^4-(8m^2+4m-2)i\xi_n^2+i}{8(\xi_n-i)^{m+4}(\xi_n+i)^{m+2}}Z_nh'(0){\rm tr}[\texttt{id}]\nonumber\\
:=&-\sum_{j,l=1}^{n-1}X_jY_l\xi_j\xi_l\frac{M_1}{8(\xi_n-i)^{m+4}(\xi_n+i)^{m+2}}Z_nh'(0){\rm tr}[\texttt{id}]\nonumber\\
&+X_nY_n\frac{M_1}{8(\xi_n-i)^{m+4}(\xi_n+i)^{m+2}}Z_nh'(0){\rm tr}[\texttt{id}].
\end{align}
Therefore, we get
\begin{align}\label{c22}
\widetilde{\Phi}_4=&-i\int_{|\xi'|=1}\int^{+\infty}_{-\infty} \sum_{j,l=1}^{n-1}X_jY_l\xi_j\xi_l\frac{M_1}{8(\xi_n-i)^{m+4}(\xi_n+i)^{m+2}}Z_nh'(0){\rm tr}[\texttt{id}]d\xi_n\sigma(\xi')dx'\nonumber\\
&+i\int_{|\xi'|=1}\int^{+\infty}_{-\infty} X_nY_n\frac{M_1}{8(\xi_n-i)^{m+4}(\xi_n+i)^{m+2}}Z_nh'(0){\rm tr}[\texttt{id}]d\xi_n\sigma(\xi')dx'\nonumber\\
=&-\frac{i}{16}Vol(S^{n-2})g(X^T,Y^T)Z_nh'(0)2^m\int_{\Gamma^{+}} \frac{M_1}{m(\xi_n-i)^{m+4}(\xi_n+i)^{m+2}}d\xi_ndx'\nonumber\\
&+\frac{i}{8}Vol(S^{n-2})X_nY_nZ_nh'(0)2^m\int_{\Gamma^{+}} \frac{M_1}{(\xi_n-i)^{m+4}(\xi_n+i)^{m+2}}d\xi_ndx'\nonumber\\
=&-\frac{i}{16}Vol(S^{n-2})g(X^T,Y^T)Z_nh'(0)2^m\frac{2\pi i}{(m+3)!}\left[\frac{M_1}{m(\xi_n+i)^{m+2}}\right]^{(m+3)}\bigg|_{\xi_n=i}dx'\nonumber\\
&+\frac{i}{8}Vol(S^{n-2})X_nY_nZ_nh'(0)2^m\frac{2\pi i}{(m+3)!}\left[\frac{M_1}{(\xi_n+i)^{m+2}}\right]^{(m+3)}\bigg|_{\xi_n=i}dx'\nonumber\\
:=&\Big(g(X^T,Y^T)\frac{\pi}{8(m+3)!}G_0-X_nY_n\frac{\pi}{4(m+3)!}G_1\Big)Vol(S^{n-2})Z_nh'(0)2^mdx',
\end{align}
where
\begin{align}
G_0=&\left[\frac{M_1}{m(\xi_n+i)^{m+2}}\right]^{(m+3)}\bigg|_{\xi_n=i}\nonumber\\
=&\frac{(32m^5+96m^4+88m^3-56m^2-22m+9)(m+1)(m+2)\cdot\cdot\cdot(2m-1)}{2^{2m+3}};\nonumber\\
G_1=&\left[\frac{M_1}{(\xi_n+i)^{m+2}}\right]^{(m+3)}\bigg|_{\xi_n=i}\nonumber\\
=&\frac{(32m^5+96m^4+88m^3-56m^2-22m+9)m(m+1)\cdot\cdot\cdot(2m-1)}{2^{2m+3}}.
\end{align}\\

\noindent {\bf  Case (b-V)}~$r=-1,~l=-(2m-1),~k=j=|\alpha|=0$.\\

\noindent By (\ref{c2}), we get
\begin{align}\label{c23}
\widetilde{\Phi}_5=-i\int_{|\xi'|=1}\int^{+\infty}_{-\infty}{\rm tr} [\pi^+_{\xi_n}\sigma_{-1}(\nabla_X^{S(TM)}\nabla_Y^{S(TM)}D^{-2})\times
\partial_{\xi_n}\sigma_{-(2m-1)}(D^{-(2m-1)})](x_0)d\xi_n\sigma(\xi')dx'.
\end{align}
By (3.27) in \cite{Wa7}, we have
\begin{equation}\label{62}
\partial_{\xi_n}\sigma_{-(2m-1)}(D^{-(2m-1)})(x_0)=i\bigg[\frac{(1-2m)\xi_n^2+1}{(1+\xi_n^2)^{m+1}}c(dx_n)-\frac{2m\xi_n}{(1+\xi_n^2)^{m+1}}c(\xi')\bigg].
\end{equation}
By Lemma \ref{lem3} and Lemma \ref{lem356}, we have
\begin{align}\label{ggg}
\sigma_{-1}(c(Z)\nabla_X^{S(TM)}\nabla_Y^{S(TM)}D^{-2})=&\sigma_{1}(c(Z)\nabla_X^{S(TM)}\nabla_Y^{S(TM)})\sigma_{-2}(D^{-2})+\sigma_{2}(c(Z)\nabla_X^{S(TM)}\nabla_Y^{S(TM)})\sigma_{-3}(D^{-2})\nonumber\\
&+\sum_{j=1}^n\partial_{\xi_j}\sigma_{2}(c(Z)\nabla_X^{S(TM)}\nabla_Y^{S(TM)})D_{x_j}[\sigma_{-2}(D^{-2})]\nonumber\\
:=&B^1+B^2+B^3,
\end{align}
where
\begin{align}
B^1(x_0)=&\sqrt{-1}c(Z)\bigg(\sum_{j,l=1}^nX_j\frac{\partial{Y_l}}{\partial_{x_j}}\xi_l+\sum_jA(Y)X_j\xi_j+\sum_lA(X)Y_l\xi_l\bigg)|\xi|^{-2};\nonumber\\
B^2(x_0)=&\sqrt{-1}c(Z)\sum_{j,l=1}^{n}X_jY_l\xi_j\xi_l\bigg(|\xi|^{-4}\xi_k(\Gamma^k-2\sigma^k)+|\xi|^{-6}2\xi^j\xi_\alpha\xi_\beta\partial_jg^{\alpha\beta}\bigg);\nonumber\\
B^3(x_0)=&\sqrt{-1}c(Z)\bigg(\sum_{l=1}^{n-1}(X_nY_l+X_lY_n)\xi_l+2X_nY_n\xi_n\bigg)h'(0)|\xi'|^2.
\end{align}
Firstly, the following results are obtained by further calculation of $B^1(x_0$)
\begin{align}
B^1(x_0)=&i\sum_{j=1}^n\sum_{l=1}^{n-1}X_j\frac{\partial {Y_l}}{\partial {x_j}}\xi_l\frac{c(Z)}{|\xi|^2}+i\sum_{j=1}^nX_j\frac{\partial {Y_n}}{\partial {x_j}}\xi_n\frac{c(Z)}{|\xi|^2}+ic(Z)\sum_{j=1}^{n-1}A(Y)X_j\xi_j\frac{1}{|\xi|^2}\nonumber\\
&+ic(Z)A(Y)X_n\xi_n\frac{1}{|\xi|^2}+ic(Z)\sum_{l=1}^{n-1}A(X)Y_l\xi_l\frac{1}{|\xi|^2}+ic(Z)A(X)Y_n\xi_n\frac{1}{|\xi|^2}.
\end{align}
If we omit some items that have no contribution for computing ${\bf \widetilde{\Phi}_5}$, we obtain
\begin{align}\label{64}
\pi^+_{\xi_n}B^1(x_0)=&\sum_{j=1}^nX_j\frac{\partial {Y_n}}{\partial {x_j}}\frac{ic(Z)}{2(\xi_n-i)}+X_n\frac{ic(Z)A(Y)}{2(\xi_n-i)}+Y_n\frac{ic(Z)A(X)}{2(\xi_n-i)}.
\end{align}
We note that $i<n,~\int_{|\xi'|=1}\{\xi_{i_{1}}\xi_{i_{2}}\cdots\xi_{i_{2q+1}}\}\sigma(\xi')=0$, and by (\ref{a49111}) , we have
\begin{equation}
{\rm tr} [\pi^+_{\xi_n}B^1\times
\partial_{\xi_n}\sigma_{-(2m-1)}(D^{-(2m-1)})](x_0)|_{|\xi'|=1}=\sum_{j=1}^nX_j\frac{\partial {Y_n}}{\partial {x_j}}\frac{(1-2m)\xi_n^2+1}{2(\xi_n-i)^{m+2}(\xi_n+i)^{m+1}}Z_n{\rm tr}[\texttt{id}].
\end{equation}
Then
\begin{align}\label{65}
&-i\int_{|\xi'|=1}\int^{+\infty}_{-\infty}{\rm tr} [\pi^+_{\xi_n}B^1\times
\partial_{\xi_n}\sigma_{-(2m-1)}(D^{-(2m-1)})](x_0)d\xi_n\sigma(\xi')dx'\nonumber\\
=&-\frac{i}{2}\int_{|\xi'|=1}\int^{+\infty}_{-\infty}\sum_{j=1}^nX_j\frac{\partial {Y_n}}{\partial {x_j}}\frac{(1-2m)\xi_n^2+1}{(\xi_n-i)^{m+2}(\xi_n+i)^{m+1}}Z_n{\rm tr}[\texttt{id}]d\xi_n\sigma(\xi')dx'\nonumber\\
=&-\frac{i}{2}Vol(S^{n-2})X(Y_n)Z_n2^m\int_{\Gamma^{+}}\frac{(1-2m)\xi_n^2+1}{(\xi_n-i)^{m+2}(\xi_n+i)^{m+1}}d\xi_ndx'\nonumber\\
=&-\frac{i}{2}Vol(S^{n-2})X(Y_n)Z_n2^m\frac{2\pi i}{(m+1)!}\left[\frac{(1-2m)\xi_n^2+1}{(\xi_n+i)^{m+1}}\right]^{(m+1)}\bigg|_{\xi_n=i}dx'\nonumber\\
:=&X(Y_n)Z_n\frac{\pi}{(m+1)!}H_0Vol(S^{n-2})2^mdx',
\end{align}
where
\begin{equation}
H_0=\left[\frac{(1-2m)\xi_n^2+1}{(\xi_n+i)^{m+1}}\right]^{(m+1)}\bigg|_{\xi_n=i}=\frac{m(m+1)\cdot\cdot\cdot(2m-1)}{2^{2m}}.
\end{equation}
Secondly, by (3.26) in \cite{Wa6}, further calculation leads to new results
\begin{align}
B^2(x_0)=&ic(Z)\sum_{j,l=1}^{n}X_jY_l\xi_j\xi_lh'(0)\bigg[\frac{\sum_{k<n}\xi_kc(e_k)c(e_n)}{2(1+\xi_n^2)^2}+\frac{m\xi_n^3+(m+2)\xi_n}{(1+\xi_n^2)^3}\bigg].
\end{align}
By the relation of the Clifford action and ${\rm tr}{(AB)}={\rm tr }{(BA)}$,  we have the equalities:
\begin{align}\label{a491}
&{\rm tr}[c(Z)\sum_{k<n}\xi_kc(e_k)c(e_n)c(\xi')]=-\sum_{k<n}\xi_k^2Z_n{\rm tr}[\texttt{id}];\nonumber\\
&{\rm tr}[c(Z)\sum_{k<n}\xi_kc(e_k)c(e_n)c(dx_n)]=g(Z,\xi'){\rm tr}[\texttt{id}].
\end{align}
If we omit some items that have no contribution for computing ${\bf \widetilde{\Phi}_5}$, we obtain
\begin{align}\label{64}
\pi^+_{\xi_n}B^2(x_0)=&\sum_{j,l=1}^{n-1}X_jY_l\xi_j\xi_lh'(0)\bigg[\frac{(\xi_n-2i)c(Z)\sum_{k<n}\xi_kc(e_k)c(e_n)}{8(\xi_n-i)^2}+\frac{[(2m+1)\xi_n-(2m+3)i]c(Z)}{8(\xi_n-i)^3}\bigg]\nonumber\\
&+X_nY_nh'(0)\bigg[\frac{\xi_nc(Z)\sum_{k<n}\xi_kc(e_k)c(e_n)}{8(\xi_n-i)^2}+\frac{[4mi\xi_n^2+(6m+3)\xi_n-(2m+1)i]c(Z)}{8(\xi_n-i)^3}\bigg].\nonumber\\
\end{align}
So
\begin{align}\label{ggg32}
&{\rm tr} [\pi^+_{\xi_n}B^2\times
\partial_{\xi_n}\sigma_{-(2m-1)}(D^{-(2m-1)})](x_0)|_{|\xi'|=1}\nonumber\\
=&\sum_{j,l=1}^{n-1}X_jY_l\xi_j\xi_l\frac{mi\xi_n^2+2m\xi_n}{8(\xi_n-i)^{m+3}(\xi_n+i)^{m+1}}\sum_{k<n}\xi_k^2Z_nh'(0){\rm tr}[\texttt{id}]\nonumber\\
&+\sum_{j,l=1}^{n-1}X_jY_l\xi_j\xi_l\frac{(4m^2-1)i\xi_n^3+(4m^2+4m-3)\xi_n^2-(2m+1)i\xi_n-(2m+3)}{8(\xi_n-i)^{m+4}(\xi_n+i)^{m+1}}Z_nh'(0){\rm tr}[\texttt{id}]\nonumber\\
&+X_nY_n\frac{mi\xi_n^2}{8(\xi_n-i)^{m+3}(\xi_n+i)^{m+1}}\sum_{k<n}\xi_k^2Z_nh'(0){\rm tr}[\texttt{id}]\nonumber\\
&+X_nY_n\frac{(-8m^2+4m)\xi_n^4+(12m^2-3)i\xi_n^3+(4m^2+4m-1)\xi_n^2-(6m+3)i\xi_n-(2m+1)}{8(\xi_n-i)^{m+4}(\xi_n+i)^{m+1}}Z_nh'(0){\rm tr}[\texttt{id}]\nonumber\\
:=&\sum_{j,l=1}^{n-1}X_jY_l\xi_j\xi_l\bigg[\sum_{k<n}\xi_k^2\frac{mi\xi_n^2+2m\xi_n}{8(\xi_n-i)^{m+3}(\xi_n+i)^{m+1}}+\frac{M_2}{8(\xi_n-i)^{m+4}(\xi_n+i)^{m+1}}\bigg]Z_nh'(0){\rm tr}[\texttt{id}]\nonumber\\
&+X_nY_n\bigg[\sum_{k<n}\xi_k^2\frac{mi\xi_n^2}{8(\xi_n-i)^{m+3}(\xi_n+i)^{m+1}}+\frac{M_3}{8(\xi_n-i)^{m+4}(\xi_n+i)^{m+1}}\bigg]Z_nh'(0){\rm tr}[\texttt{id}].
\end{align}
Then
\begin{align}\label{65}
&-i\int_{|\xi'|=1}\int^{+\infty}_{-\infty}{\rm tr} [\pi^+_{\xi_n}B^2\times
\partial_{\xi_n}\sigma_{-(2m-1)}(D^{-(2m-1)})](x_0)d\xi_n\sigma(\xi')dx'\nonumber\\
=&-i\int_{|\xi'|=1}\int^{+\infty}_{-\infty}\sum_{j,l=1}^{n-1}X_jY_l\xi_j\xi_l\sum_{k<n}\xi_k^2\frac{mi\xi_n^2+2m\xi_n}{8(\xi_n-i)^{m+3}(\xi_n+i)^{m+1}}Z_nh'(0){\rm tr}[\texttt{id}]d\xi_n\sigma(\xi')dx'\nonumber\\
&-i\int_{|\xi'|=1}\int^{+\infty}_{-\infty}\sum_{j,l=1}^{n-1}X_jY_l\xi_j\xi_l\frac{M_2}{8(\xi_n-i)^{m+4}(\xi_n+i)^{m+1}}Z_nh'(0){\rm tr}[\texttt{id}]d\xi_n\sigma(\xi')dx'\nonumber\\
&-i\int_{|\xi'|=1}\int^{+\infty}_{-\infty}X_nY_n\sum_{k<n}\xi_k^2\frac{mi\xi_n^2}{8(\xi_n-i)^{m+3}(\xi_n+i)^{m+1}}Z_nh'(0){\rm tr}[\texttt{id}]d\xi_n\sigma(\xi')dx'\nonumber\\
&-i\int_{|\xi'|=1}\int^{+\infty}_{-\infty}X_nY_n\frac{M_3}{8(\xi_n-i)^{m+4}(\xi_n+i)^{m+1}}Z_nh'(0){\rm tr}[\texttt{id}]d\xi_n\sigma(\xi')dx'\nonumber\\
=&-\frac{3i}{32}Vol(S^{n-2})g(X^T,Y^T)Z_nh'(0)2^{m}\int_{\Gamma^{+}}\frac{i\xi_n^2+2\xi_n}{(m+1)(\xi_n-i)^{m+3}(\xi_n+i)^{m+1}}d\xi_ndx'\nonumber\\
&-\frac{i}{16}Vol(S^{n-2})g(X^T,Y^T)Z_nh'(0)2^{m}\int_{\Gamma^{+}}\frac{M_2}{m(\xi_n-i)^{m+4}(\xi_n+i)^{m+1}}d\xi_ndx'\nonumber\\
&-\frac{i}{16}Vol(S^{n-2})X_nY_nZ_nh'(0)2^{m}\int_{\Gamma^{+}}\frac{i\xi_n^2}{(\xi_n-i)^{m+3}(\xi_n+i)^{m+1}}d\xi_ndx'\nonumber\\
&-\frac{i}{8}Vol(S^{n-2})X_nY_nZ_nh'(0)2^{m}\int_{\Gamma^{+}}\frac{M_3}{(\xi_n-i)^{m+4}(\xi_n+i)^{m+1}}d\xi_ndx'\nonumber\\
=&-\frac{3i}{32}Vol(S^{n-2})g(X^T,Y^T)Z_nh'(0)2^{m}\frac{2\pi i}{(m+2)!}\left[\frac{i\xi_n^2+2\xi_n}{(m+1)(\xi_n+i)^{m+1}}\right]^{(m+2)}\bigg|_{\xi_n=i}dx'\nonumber\\
&-\frac{i}{16}Vol(S^{n-2})g(X^T,Y^T)Z_nh'(0)2^{m}\frac{2\pi i}{(m+3)!}\left[\frac{M_2}{m(\xi_n+i)^{m+1}}\right]^{(m+3)}\bigg|_{\xi_n=i}dx'\nonumber\\
&-\frac{i}{16}Vol(S^{n-2})X_nY_nZ_nh'(0)2^{m}\frac{2\pi i}{(m+2)!}\left[\frac{i\xi_n^2}{(\xi_n+i)^{m+1}}\right]^{(m+2)}\bigg|_{\xi_n=i}dx'\nonumber\\
&-\frac{i}{8}Vol(S^{n-2})X_nY_nZ_nh'(0)2^{m}\frac{2\pi i}{(m+3)!}\left[\frac{M_3}{(\xi_n+i)^{m+1}}\right]^{(m+3)}\bigg|_{\xi_n=i}dx'\nonumber\\
:=&\bigg\{g(X^T,Y^T)\Big(\frac{3\pi}{16(m+2)!}H_1+\frac{\pi}{8(m+3)!}H_2\Big)+X_nY_n\Big(\frac{\pi}{8(m+2)!}H_3\nonumber\\
&+\frac{\pi}{4(m+3)!}H_4\Big)\bigg\}Vol(S^{n-2})Z_nh'(0)2^{m}dx',
\end{align}
where
\begin{align}
H_1=&\left[\frac{i\xi_n^2+2\xi_n}{(m+1)(\xi_n+i)^{m+1}}\right]^{(m+2)}\bigg|_{\xi_n=i}=-\frac{(5m+1)(m+2)(m+3)\cdot\cdot\cdot(2m)}{2^{2m+2}};\nonumber\\
H_2=&\left[\frac{M_2}{m(\xi_n+i)^{m+1}}\right]^{(m+3)}\bigg|_{\xi_n=i}=\frac{3(2m^2+11m+5)m(m+1)\cdot\cdot\cdot(2m-1)}{2^{2m}};\nonumber\\
H_3=&\left[\frac{i\xi_n^2}{(\xi_n+i)^{m+1}}\right]^{(m+2)}\bigg|_{\xi_n=i}=-\frac{(m-1)m\cdot\cdot\cdot(2m-1)}{2^{2m+1}};\nonumber\\
H_4=&\left[\frac{M_3}{(\xi_n+i)^{m+1}}\right]^{(m+3)}\bigg|_{\xi_n=i}=-\frac{(2m^3+m^2+45m)m(m+1)\cdot\cdot\cdot(2m-1)}{2^{2m}}.
\end{align}
Thirdly, for $B^3$, we get
\begin{align}\label{66pp}
\pi^+_{\xi_n} B^3(x_0)=2iX_nY_nc(Z)h'(0)|\xi'|^2.
\end{align}
Moreover
\begin{align}
{\rm tr} [\pi^+_{\xi_n}B^3\times
\partial_{\xi_n}\sigma_{-(2m-1)}(D^{-(2m-1)})](x_0)|_{|\xi'|=1}=X_nY_n\frac{(2-4m)\xi_n^2+2}{(\xi_n-i)^{m+1}(\xi_n+i)^{m+1}}Z_nh'(0){\rm tr}[\texttt{id}].
\end{align}
Then
\begin{align}
&-i\int_{|\xi'|=1}\int^{+\infty}_{-\infty}{\rm tr} [\pi^+_{\xi_n}B^3\times
\partial_{\xi_n}\sigma_{-(2m-1)}(D^{-(2m-1)})](x_0)d\xi_n\sigma(\xi')dx'\nonumber\\
=&-i\int_{|\xi'|=1}\int^{+\infty}_{-\infty}X_nY_n\frac{(2-4m)\xi_n^2+2}{(\xi_n-i)^{m+1}(\xi_n+i)^{m+1}}Z_nh'(0){\rm tr}[\texttt{id}]d\xi_n\sigma(\xi')dx'\nonumber\\
=&-iVol(S^{n-2})X_nY_nZ_nh'(0)2^m\int_{\Gamma^+}\frac{(2-4m)\xi_n^2+2}{(\xi_n-i)^{m+1}(\xi_n+i)^{m+1}}d\xi_ndx'\nonumber\\
=&-iVol(S^{n-2})X_nY_nZ_nh'(0)2^m\frac{2\pi i}{m!}\left[\frac{(2-4m)\xi_n^2+2}{(\xi_n+i)^{m+1}}\right]^{(m)}\bigg|_{\xi_n=i}dx'\nonumber\\
:=&X_nY_n\frac{2\pi}{m!}H_5Vol(S^{n-2})Z_nh'(0)2^mdx',
\end{align}
where
\begin{align}
H_5=\left[\frac{(2-4m)\xi_n^2+2}{(\xi_n+i)^{m+1}}\right]^{(m+2)}\bigg|_{\xi_n=i}=0.
\end{align}
Therefore, we obtain
\begin{align}\label{6666}
\widetilde{\Phi}_5=&-i\int_{|\xi'|=1}\int^{+\infty}_{-\infty}{\rm tr} [\pi^+_{\xi_n}(B^1+B^2+B^3)\times
\partial_{\xi_n}\sigma_{-(2m-1)}(D^{-(2m-1)})](x_0)d\xi_n\sigma(\xi')dx'\nonumber\\
=&\bigg\{X(Y_n)Z_n\frac{\pi}{(m+1)!}H_0+g(X^T,Y^T)Z_nh'(0)\Big(\frac{3\pi}{16(m+2)!}H_1+\frac{\pi}{8(m+3)!}H_2\Big)\nonumber\\
&+X_nY_nZ_nh'(0)\Big(\frac{\pi}{8(m+2)!}H_3+\frac{\pi}{4(m+3)!}H_4\Big)\bigg\}Vol(S^{n-2})2^mdx'.
\end{align}
Now $\widetilde{\Phi}$ is the sum of the {\bf  Case (b-I)}-{\bf Case (b-V)}. Therefore, we get
\begin{align}\label{c28}
\widetilde{\Phi}&=\sum_{i=1}^5\widetilde{\Phi}_i\nonumber\\
=&\bigg\{\partial_{x_n}\Big(g(X^T,Y^T)Z_n\Big)\frac{-\pi i}{2(m+2)!}E_0+\partial_{x_n}\Big((X_nY_n)Z_n\Big)\frac{\pi i}{(m+2)!}E_1\nonumber\\
&+g(X^T,Y^T)Z_nh'(0)\Big(\frac{\pi i}{4(m+3)!}E_2-\frac{\pi}{2(m+3)!}F_0+\frac{\pi}{8(m+3)!}G_0+\frac{3\pi}{16(m+2)!}H_1+\frac{\pi}{8(m+3)!}H_2\Big)\nonumber\\
&+X_nY_nZ_nh'(0)\Big(\frac{\pi i}{2(m+3)!}E_3+\frac{\pi i}{(m+3)!}F_1-\frac{\pi}{4(m+3)!}G_1+\frac{\pi}{8(m+2)!}H_3+\frac{\pi}{4(m+3)!}H_4\Big)\nonumber\\
&+X(Y_n)Z_n\frac{\pi}{(m+1)!}H_0\bigg\}Vol(S^{n-2})2^mdx'.
\end{align}
Then, by (\ref{c121}) and (\ref{c28}), we obtain the following theorem.
\begin{thm}\label{cthmb1}
Let $M$ be an $n=(2m+1)$-dimensional oriented
compact spin manifold with the boundary $\partial M$, then we get the following equality:
\begin{align}
\label{cb2773}
&\widetilde{{\rm Wres}}[\pi^+(c(Z)\nabla_X^{S(TM)}\nabla_Y^{S(TM)}D^{-2})\circ\pi^+(D^{-(2m-1)})]\nonumber\\
=&\int_{\partial M}\bigg\{\partial_{x_n}\Big(g(X^T,Y^T)Z_n\Big)\frac{-\pi i}{2(m+2)!}E_0+\partial_{x_n}\Big((X_nY_n)Z_n\Big)\frac{\pi i}{(m+2)!}E_1\nonumber\\
&+g(X^T,Y^T)Z_nh'(0)\Big(\frac{\pi i}{4(m+3)!}E_2-\frac{\pi}{2(m+3)!}F_0+\frac{\pi}{8(m+3)!}G_0+\frac{3\pi}{16(m+2)!}H_1+\frac{\pi}{8(m+3)!}H_2\Big)\nonumber\\
&+X_nY_nZ_nh'(0)\Big(\frac{\pi i}{2(m+3)!}E_3+\frac{\pi i}{(m+3)!}F_1-\frac{\pi}{4(m+3)!}G_1+\frac{\pi}{8(m+2)!}H_3+\frac{\pi}{4(m+3)!}H_4\Big)\nonumber\\
&+X(Y_n)Z_n\frac{\pi}{(m+1)!}H_0\bigg\}Vol(S^{n-2})2^md{\rm Vol_{M}}.
\end{align}\\
\end{thm}

\noindent {\small\textbf{Acknowledgements}} This work was supported by NSFC. 11771070. The authors thank the referee for his (or her) careful reading and helpful comments.\\

\noindent {\small\textbf{Data availability statement}} The authors confirm that the data supporting the findings of this study
are available within the article.

\section*{Declarations}
\noindent {\small\textbf{Conflict of interest}} The authors state that there is no conflict of interest.

\end{document}